%% file: main.tex
\documentclass[12pt, a4paper]{article}

\usepackage{lscape}

\usepackage[version=4]{mhchem}
\usepackage{amsmath,amsfonts,amssymb}
\usepackage{eurosym}
\usepackage[T1]{fontenc}
\usepackage[latin1]{inputenc}
\usepackage{setspace}
\usepackage{xcolor}
\usepackage{textcomp}
\usepackage{gensymb}
\usepackage{layouts}
\usepackage{mathrsfs}

\usepackage{tabularx}
\usepackage[flushleft]{threeparttable}
\usepackage{makecell}
\usepackage{multirow}
\usepackage{exscale}
\usepackage{array}
\usepackage{dcolumn}
\usepackage{booktabs}
\usepackage{csvsimple}

\newcolumntype{L}[1]{>{\raggedright\arraybackslash}p{#1}}

\usepackage{tikz}
\usetikzlibrary{positioning, calc}
\usetikzlibrary[backgrounds]
\usetikzlibrary{shapes.geometric, arrows}
\usetikzlibrary{arrows.meta, decorations.markings}
\usetikzlibrary{external}
\usepackage{pgfplots}
\usepackage{pgfplotstable}
\usepgfplotslibrary{groupplots}
\pgfplotsset{compat=1.7}
\usepackage[skip=0.3\baselineskip]{caption}
\usepackage{subcaption}
\usepackage{graphicx} 		
\usepackage{svg}
\usepackage{wrapfig}
\usepackage{placeins}

\usepackage[numbers]{natbib}
\usepackage{ragged2e} 						
\usepackage{float}	
\usepackage{acronym}
\usepackage{hyperref}
\usepackage{etoolbox}
\usepackage{needspace}
\usepackage{siunitx}
\usepackage{rotating}
\usepackage{import}
\usepackage{ifthen}
\usepackage[nottoc,numbib]{tocbibind}
\usepackage{enumitem}
\usepackage{comment}

\AtBeginEnvironment{tabular}{\scriptsize}
\AtBeginEnvironment{tabularx}{\scriptsize}
\DeclareMathSizes{10}{10}{4}{2}

\usepackage[margin=1in]{geometry}
 
\newcolumntype{d}{D{.}{.}{-1}}
\newcolumntype{Y}{>{\centering\arraybackslash}X}
\newcolumntype{P}[1]{>{\centering\arraybackslash}p{#1}}

\definecolor{myblue}{RGB}{14, 48, 96}
\definecolor{myred}{RGB}{166, 38, 35}
\definecolor{myyellow}{RGB}{209, 145, 62}
\definecolor{mylightblue}{RGB}{19, 83, 144}

\hypersetup{	
	linkcolor=myblue,
	colorlinks = true,
	citecolor=myblue
}

\makeatletter
\newif\if@in@acrolist
\AtBeginEnvironment{acronym}{\@in@acrolisttrue}
\newrobustcmd{\LU}[2]{\if@in@acrolist#1\else#2\fi}

\newcommand{\ACF}[1]{{\@in@acrolisttrue\acf{#1}}}
\makeatother

\newcommand*{\SmallIndent}{\hspace*{0.2cm}}%

\title{Unleashing the full potential of the North Sea - Identifying key energy infrastructure synergies for 2030 and 2040}

\author{Jan Wiegner, Madeleine Gibescu, Matteo Gazzani}

\begin{document}
\maketitle

\textbf{Abstract}

Policy efforts have primarily focused on expanding variable renewable energy sources (vRES) to meet ambitious carbon emission reduction targets. The integration of high shares of renewables into existing and future energy systems is central to both policy making and research, focusing on the need for balancing options between variable renewable energy sources and demand. In this work we analyze and compare three key integration measures: grid expansions, electricity storage, and the role of production, storage and transport of low-carbon hydrogen. We focus on their potential to reduce emissions and energy system costs,  individually and in combination with each other. This allows for exploring synergies between these three integration measures. Additionally, we study the consequences of regulatory constraints on their implementation. We take the North Sea as an exemplary region with ambitious 2030-2040 targets for offshore wind developments. The projections on installed generation and grid capacities, along with demand estimates from the Ten Year Network Development Plan (TYNDP) 2022 by the European Network of Transmission System Operators, serve as a starting point for our energy system model. This starting model can then be further expanded, in a cost- and  emission- minimization fashion, with the three integration measures. Our findings show that electricity grid expansions across the North Sea are a no-regret measure lowering costs, emissions and required renewable expansions in the region. The production of hydrogen and its direct use in industry has a lower cost reduction potential and emission reduction potential, while hydrogen storage and transport have little to no additional value. In the short term (2030), electricity storage can help to reduce emissions, but it is not cost competitive. In the longer term (2040), storage can help to balance investments in vRES assets by providing additional flexibility to the system. Combining the three integration measures provides additional benefits. The highest emission reductions can be achieved by combining electricity storage with an expansion of the grid. The highest economic benefits can be achieved with a combination of grid expansions and hydrogen production for direct use in industry. 

\newpage

\newpage

%
\newpage
\thispagestyle{empty}
\section*{List of Abbreviations}
\begin{acronym}[Bash]
    \acro{CCS}{Carbon Capture \& Storage}
    \acro{ENTSO-E}{European Network of Transmission System Operators for Electricity}
    \acro{ERAA}{European Resource Adequacy Assessment}
    \acro{HVDC}{High Voltage Direct Current}
    \acro{TSO}{Transmission System Operator}
    \acro{TYNDP}{Ten Year Network Development Plan}
    \acro{vRES}{Variable Renewable Energy Sources}
\end{acronym}
\section*{List of Mathematical Symbols}
\begin{table}[H]
    \scriptsize
    \begin{tabularx}{\columnwidth}{lX}
    \midrule
    \textbf{Indices and Sets} & \\
        $\mathrm{n \in N}$ & Set of nodes\\
        $\mathrm{r \in R = \{\mathrm{el, h2, ng}}\}$ & Set of energy carriers\\
        $\mathrm{t \in T}$ & Set of timesteps\\
        $\mathrm{i \in I_{n}}$ & Set of technologies at node $n$\\
        $\mathrm{g \in G}$ & Set of networks\\
        $\mathrm{l \in L_{g}}$ & Set of branches of network $g$\\
    \midrule
    \textbf{Subscripts}\\
        $\mathrm{exp}$ & Export\\
        $\mathrm{imp}$ & Import\\
        $\mathrm{out}$ & output of a technology or outflow from a node to a line of a network\\
        $\mathrm{in}$ & input to a technology or inflow to a node from a line of a network\\
        $\mathrm{c}$ & Up front investment cost\\
        $\mathrm{f}$ & Fixed cost\\
        $\mathrm{v}$ & Variable cost\\
        $\mathrm{CO2}$ & Carbon emissions\\
        $\mathrm{cons}$ & Energy consumption of networks\\
        $\mathrm{netw}$ & Networks\\
        $\mathrm{tech}$ & Technologies\\
    \midrule
    \textbf{Variables} & \\
        $E$ & Emissions\\
        $X$ & Technology input or output/ Import or Export\\
        $S$ & Technology or network branch size\\
        $F$ & Network flow\\
        $C$ & Cost\\
        $R$ & Revenues\\
    \midrule
    \textbf{Parameters} & \\
        $\mathrm{c}$ & Cost parameter\\
        $\mathrm{e}$ & Emission factor\\
        $\mathrm{D}$ & Demand\\
        $\mathrm{d}$ & Distance between two nodes\\
        $\mathrm{\gamma}$ & Cost parameter for branch costs of a network\\
    \midrule
    \end{tabularx}
\end{table}
\newpage

\setcounter{page}{1}
\section{Introduction}
\acf{vRES} such as wind and solar power are increasingly contributing to global electricity generation. However, their intermittent nature presents significant challenges when it comes to balancing supply and demand in an electricity system. The temporal mismatch between \ac{vRES} supply and electricity demand may lead to curtailment of carbon-free electricity at times of oversupply, but also requires dispatchable back-up capacities during shortages that are often based on fossil fuels. So far, policymakers have extensively focused on defining ambitious goals for renewable capacity expansions to reduce the carbon footprint of the energy system. The required \ac{vRES} capacities for deep decarbonization goals of the whole economy (50\% reduction by 2030 and net zero emissions by 2050) must not only meet current demand but also accommodate the significant additional demand arising from the anticipated electrification of heating, transport, and industry sectors \cite{EuropeanCommission2019}. While progress toward reaching the targets for capacity expansions of \ac{vRES} are a clear focus of policy making, the integration of high shares of \ac{vRES} into the existing energy system has only recently begun to gain momentum, driven by the already visible increase in congestion and real-time imbalances in electricity grids. As a result of support schemes and regulations mainly targeting the expansion of \ac{vRES} capacities, high-share \ac{vRES} systems in Europe have reached or will soon reach their integration limits \cite{EUAgencyfortheCooperationofEnergyRegulators2023, Davi-Arderius2024}.

In this work, we focus on technology measures for a successful integration of high shares of \ac{vRES} towards 2030 and 2040 in the North Sea region. Therefore, we take the development of demand and the expansion of generation capacities as exogenous variables -- as provided by the ENTSO-e Ten Year Network Development Plan (TYNDP) -- and focus on the desirable 'integration measures', that we will call 'supportive infrastructure'. The supportive infrastructure falls broadly into three categories: (i) electricity grid expansions ("\textit{grid expansion}"), (ii) "\textit{electricity storage}" and (iii) the production, storage, transport of hydrogen as well as its re-conversion into electricity and its direct use in industry ("\textit{hydrogen technologies}"). The facilitating role of the supportive infrastructure is paramount and has been stressed in multiple studies, finding significant potential to lower both costs and emissions of a future energy system \cite{Neumann2023, Sgobbi2016, Pickering2022, Seck2022, Martinez-Gordon2022a, Golombek2022}. However, a systematic comparison between different integration measures is missing.

In this work, we thus compare these three supportive infrastructure measures in a holistic matter, i.e. from a global societal welfare point of view and along their emission reduction potentials, when applied separately and in combinations. We focus on the North Sea region, an exemplary geographic area expecting a substantial increase in \ac{vRES} capacities by 2030, with projections of 400 GW renewable capacities from solar, onshore, and offshore wind. Additionally, the neighboring countries around the North Sea are key economic hubs, contributing significantly to Europe's economic output. This makes the North Sea region a particularly interesting area for studying the integration of \ac{vRES}. The North Sea basin is thereby gaining specifically increasing interest not only in the expansion of offshore wind capacities but also in expanding other energy infrastructure, such as hydrogen technologies, energy storage and electricity transmission infrastructure. 

In this work, we look at these three integration measures at two different points in time: (i) the short term (year 2030), assuming exogenous electricity and hydrogen demand as well as fixed conventional and renewable generation capacities as provided by TYNDP 2022 and (ii) the medium term (year 2040) with exogenous electricity and hydrogen demand and a possibility for expansion of on- and offshore \ac{vRES} capacities. We thus pose the following overarching research questions:
\begin{itemize}
    \item How can supportive infrastructure (grid expansions, electricity storage and hydrogen technologies) contribute to welfare gains and emission reduction in the short (2030) and medium term (2040)?
    \item Are there synergies between the three integration measures?
    \item What is the societal welfare effect and the impact on the emission reduction potential if a measure is infeasible due to political, legal, social or technical constraints?
\end{itemize}

Hereafter, we introduce previous work on the three integration measures and pose more specific research questions, which are discussed in the Results section.

\vspace{0.5cm}
\textbf{The role of electricity grids.}
Previous studies on the role of the electricity grid for the integration of \ac{vRES} have generally found that grid expansions not only reduce curtailment and thus emissions, but also result in high overall welfare benefits \cite{GorensteinDedecca2018, Huertas-Hernando2010, GorensteinDedecca2016, Konstantelos2017, Kristiansen2018, Houghton2016, Koivisto2020a, Neumann2023, Martinez-Gordon2022a}. This is remarkably the case for the North Sea region, where submarine cables could interconnect different countries as well as various offshore wind farms. Compared to radial connections of offshore wind farms to a national landing point, the interconnection of different wind farms and multiple landing points enables the offshore grid to serve as both a means of transport from a wind farm to multiple onshore load centers as well as to facilitate the exchange of electricity between two countries.  Gorenstein-Dedeca et al. (2018) and Hadush et al. (2014), however, stress that the welfare benefit could be highly unequally shared, resulting in the potential opposition from parties or countries with a net welfare loss or only minor benefits \cite{GorensteinDedecca2018, Hadush2014}. Neumann et al. (2023) further compare the roll-out of an integrated electricity grid with a hydrogen network and find that the expansion of the electricity grid has higher economic and environmental benefits than the investment in an interconnected hydrogen network. The highest benefit for a carbon-neutral energy system towards 2050, however, results from the co-existence of both networks \cite{Neumann2023}. To summarize, while the existing body of literature agrees that the expansion of electricity grids are economic for the successful integration of \ac{vRES}, it falls short on two dimensions: (i) the comparison to other integration measures (see overarching research questions above) and (ii) evaluating the interplay between onshore and offshore international grid expansions. Our work sheds light on these two aspects by providing possible answers to the following questions:
\begin{itemize}
    \item How important is the expansion of the offshore electricity grid compared to its onshore counterpart?
    \item How important are cross-border electricity interconnectors to realize the full potential of grid expansions?
\end{itemize}


\textbf{The role of the hydrogen technologies.}
Multiple studies of the European energy system have confirmed the future role of hydrogen in transitioning towards a carbon neutral economy \cite{Sgobbi2016, Pickering2022,Seck2022,Martinez-Gordon2022a}. 
As a result, the optimal hydrogen supply infrastructure has been studied extensively, accounting for storage, transport and production options. While low-carbon hydrogen can come from multiple sources (biomass, natural gas with CCS, nuclear electricity, or \ac{vRES} electricity), the role of hydrogen from electrolysis has been stressed frequently as a means to lower curtailment from \ac{vRES} or to increase the capacity factors of nuclear power plants \cite{Ganter2024, Kakoulaki2021, Won2017, Weimann2021, Neumann2023, Reuss2019}. Furthermore, with increasing offshore wind capacities, hydrogen production offshore has gained further interest recently. The main rationale for moving its production offshore is (i) to avoid grid congestion offshore, (ii) to enable energy storage in the form of hydrogen possibly in depleted oil and gas field and thus reusing existing infrastructure offshore, and (iii) to save on the electricity connection of new offshore wind farms by transporting hydrogen in existing pipelines to shore \cite{Wiegner2024a}. McDonagh et al. (2020) and Yan (2021) compare hydrogen production offshore to the option of feeding electricity generated from wind farms directly into the grid \cite{McDonagh2020, Dinh2021}. They find that electricity transmission is the preferred option, unless the price of green hydrogen is significantly higher then the production cost of grey hydrogen. Other studies focusing solely on hydrogen production have indicated that for hydrogen to become a competitive business model, its selling price will need to rise substantially from 2024 levels \cite{Baldi2022, Franco2021, Hou2017c}. Other studies, however, find that offshore hydrogen production is more economic than onshore production, when assuming similar investment costs for onshore and offshore electrolysis \cite{Singlitico2021b, Jang2022}. The existing body of literature falls short on a system-wide perspective of hydrogen production from offshore wind and as such, we pose the following questions in this work:
\begin{itemize}
    \item Does hydrogen production offshore provide environmental and welfare benefits comparable to its onshore production?
    \item Can hydrogen production, storage and long-distance transport as well as its re-conversion into electricity or its direct use in industry to replace blue hydrogen contribute to reduce system costs and/or to decarbonize the energy system?
\end{itemize}


\textbf{The role of electricity storage.}
Previous work on the role of electricity storage for carbon emission reduction suggest that storage is essential to reach carbon-neutrality, but it is also comparatively expensive \cite{Victoria2019, Pickering2022, Cebulla2017, Golombek2022}. Victoria et al. (2019) study the role of sector coupling (electricity, transport and heating) and electricity storage and find that sector coupling enables more economic decarbonization then electricity storage \cite{Victoria2019}. Similarly, Golombek et al. (2022) find that transmission capacities should be expanded before investing in electricity storage on the pathway towards carbon neutrality \cite{Golombek2022}. Studies on the role of offshore storage typically look at an individual wind farm in conjunction with a storage technology and find that the considered storage technologies are too expensive to be implemented \cite{Li2015, Simpson2021, Jafari2020, Ding2018, Denholm2019, Diaz-Gonzalez2012, Zhao2015, Loisel2012}. While studies typically look through an emission reduction or an economic perspective, it is important to note that electricity storage can have benefits beyond these dimensions. Parzen et al. (2022) stress that a certain storage technology might appear to be preferred from an economic perspective, but it might not be the technology with the highest system benefit (e.g. in terms of power system services)\cite{Parzen2022}. In this work, we study electricity storage from the perspective of societal welfare increases and in its role as an emission reduction measure. Thereby, we emphasize the offshore environment and compare offshore storage in proximity to wind farms to its placement onshore. As such we pose the following questions:
\begin{itemize}
    \item To what extent can electricity storage contribute to an economic integration of large-scale \ac{vRES}?
    \item What is the role of offshore electricity storage and how does it compare to onshore storage in terms of welfare and emission reduction potential?
\end{itemize}

The remainder of this paper is structured as follows: Section 2 describes the baseline system of 2030 based on the projections by \ac{ENTSO-E}. Section 3 provides a brief description of the mathematical formulation of the energy system model and the associated input data, and also defines the scenarios used hereafter. Section 4 presents the simulation results for all scenarios (2030 and 2040) and for all integration measures. Section 5 discusses the results and their limitations, followed by the conclusion in Section 6.

\FloatBarrier\newpage
\section{System Description}

\begin{figure}[ht]
    \centering
    \includegraphics[width=0.98\textwidth]{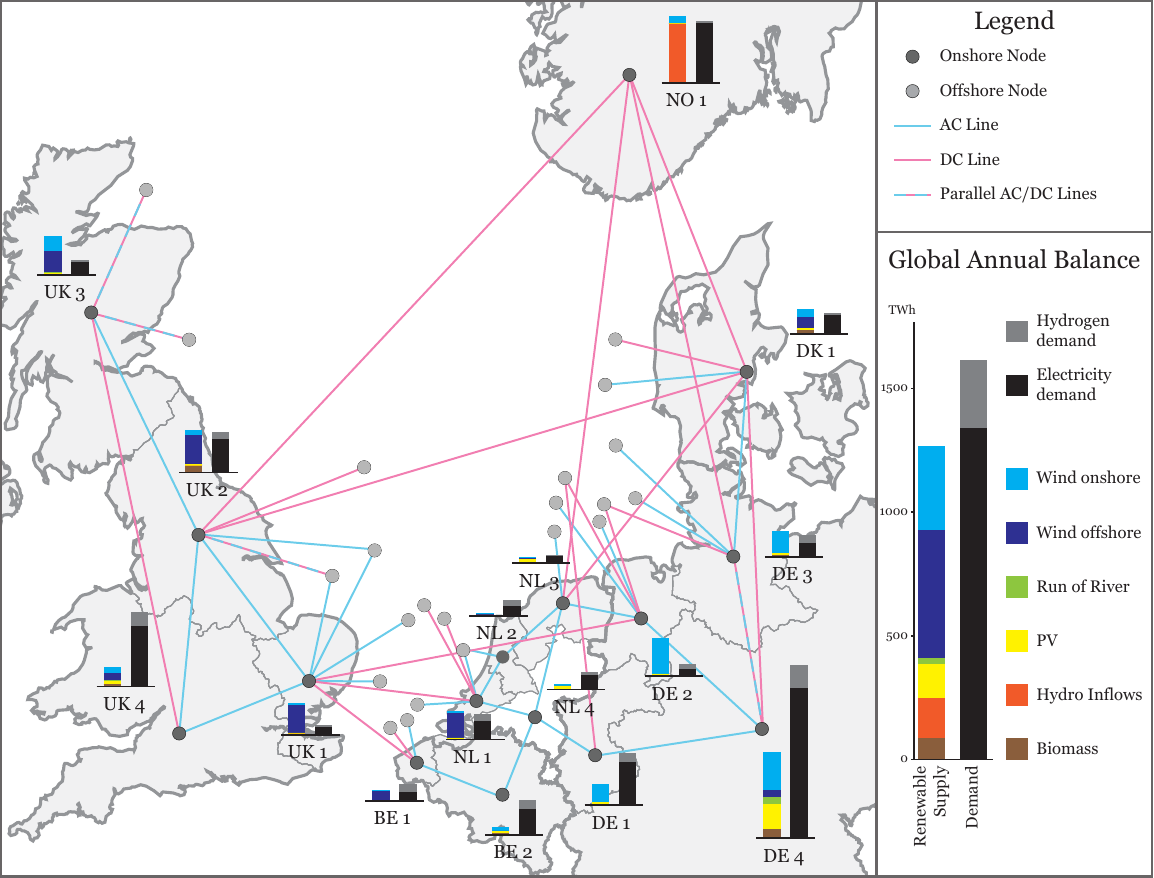}
    \caption{Illustration of the topology of the starting energy system. The vertical bars show the annual sum of theoretical supply of renewable resources without curtailment in comparison to annual demand at the same node as projected for the year 2030. Hydro inflows refers to natural water inflows into the upper reservoir of pumped hydro storage plants.}
    \label{fig:SystemTopology}
\end{figure}

We model the energy system of the North Sea neighboring states (United Kingdom, Norway, Denmark, Belgium, The Netherlands and Germany) with a focus on the electricity system in the year 2030. The countries are further divided into a number of onshore and offshore nodes to capture grid congestion problems of the onshore system. Each node is subject to an electricity demand that can be met by a number of renewable, conventional and hydro storage technologies as well as by exchange through the grid with neighboring nodes. Additionally, each node has a hydrogen demand that is supplied by blue hydrogen from steam methane reforming with carbon capture.

Energy generation and storage technologies, as well as electricity grids are modeled at their capacities as projected in 2030 by \ac{ENTSO-E} in their National Trends scenario. As such, the input data to the model already includes expected expansions of the grid as well as generation capacities. In the \textit{Reference} scenario, we optimize the operation of the reference system with no further capacity or grid additions. In all other scenarios (described in section \ref{sec:ScenarioDefinition}), the role of three different integration measures are studied: (i) additional grid expansions, (i) additional energy storage, and (iii) green hydrogen generation and storage. In these scenarios, we co-optimize the system design and operation taking into account existing infrastructure (brownfield approach). The model is formulated with an hourly resolution to capture fluctuations in renewable generation and respective storage requirements. Figure \ref{fig:SystemTopology} shows an illustration of the electricity grid topology and the demand-supply balance of the starting system. Table \ref{tab:00_InstalledCapacities} depicts the capacities of storage and generation technologies as well as networks considered.

\input{Tables/00_InstalledCapacities_Countries.tex}

\FloatBarrier\newpage
\section{Methodology}

\subsection{Energy System Model}
\subsubsection{Problem Formulation}
The energy system model used in this work is formulated as a mixed-integer linear program. It can be described in its most general form as:

\begin{equation}
\begin{split}
    \mathrm{min}_{x,y} \qquad  & f_{x,y} \\
    \mathrm{s.t.} \qquad  & \mathrm{A}x \leq \mathrm{b} \\
                          & \mathrm{C}y \leq \mathrm{d} \\
                          & \mathrm{A} \in \mathbb{R}^{m \times n}, \mathrm{C} \in \mathbb{R}^{p \times q} \\
                          & \mathrm{b} \in \mathbb{R}^{m}, \mathrm{d} \in \mathbb{R}^{p} \\
                          & x \in \mathbb{Z}^{n}, y \in \mathbb{R}^{q}
\end{split}
\end{equation}

Where $f$ is the objective function to be minimized subject to a number of inequality and equality constraints. $x$ depict all integer variables, and $y$ all real variables. Subsequently, we consider three types of optimizations: (i) emission minimization, (ii) cost minimization and (ii) cost minimization at an emission reduction target. All three types share the same constraints with regards to the energy balance, the technology performance and the network performance.

\subsubsection{Objective Function}
\textbf{Emission minimization}. The total annual emissions are composed of emissions from technologies $E_\mathrm{tec}$ and emissions from electricity imports $E_\mathrm{imp}$. In the emission minimization case, the objective function is thus:

\begin{equation}
    \label{eq:MinEmissions}
    f = E_\mathrm{tec} + E_\mathrm{imp}
\end{equation}

Where emissions from conversion and storage technologies are calculated as the sum of emissions from all nodes $\mathrm{N}$, technologies $\mathrm{I_{n}}$, carriers $\mathrm{R}$ and time slices $\mathrm{T}$. $\mathrm{e_{i,r}}$ hereby refers to the emission factor of the respective input $X_\mathrm{{in, i,r,t}}$ or output $X_\mathrm{{i,r,t}}$ to the technology $X_\mathrm{{out, i,r,t}}$.

\begin{equation}
    \label{eq:Emissions_Tec}
    E_\mathrm{tec} = \sum_{\mathrm{n \in N}} \sum_{\mathrm{i \in I_{n}}} \sum_{\mathrm{r \in R}} \sum_{\mathrm{t \in T}}  \left( \mathrm{e_{in,i,r}} X_\mathrm{{in,i,r,t}} + \mathrm{e_{out,i,r}} X_\mathrm{{out,i,r,t}}\right)
\end{equation}

Emissions from electricity imports from beyond the system boundaries are calculated as:

\begin{equation}
    \label{eq:Emissions_Imp}
    E_\mathrm{imp} = \sum_{\mathrm{n \in N}} \sum_{\mathrm{t \in T}} \left( \mathrm{e_{imp,el}}  X_\mathrm{{imp,el,n,t}} \right)
\end{equation}

\textbf{Cost minimization}. In the cost minimization problem, total annualized costs of the energy system are minimized. The total annualized cost is composed of costs for conversion and storage technologies $C_\mathrm{tec}$, network costs $C_\mathrm{netw}$, cost for imports $C_\mathrm{imp}$, and carbon costs $C_\mathrm{CO2}$. As such, the objective function is:

\begin{equation}
    \label{eq:MinCost}
    f = C_\mathrm{tec} + C_\mathrm{netw} + C_\mathrm{imp} + C_\mathrm{CO2}
\end{equation}

Costs for storage and conversion technologies $C_\mathrm{tec}$ are composed of size dependent investment cost $\mathrm{c_{c,i}} S_\mathrm{{i}}$, fixed costs calculated as a fraction of investment costs ($\mathrm{c_{f,i}} C_{c,i}$) and variable O\&M costs depending on the technology output (eq. \ref{eq:Cost_Techs}). Note that for technologies included in the starting system do not have investment costs (i.e. $\mathrm{c_{c,i}} = \mathrm{c_{f,i}} = C_{c,i} = 0$).
\begin{equation}
    \label{eq:Cost_Techs}
    C_\mathrm{tec} = \sum_{\mathrm{\mathrm{n \in N}}} \sum_{\mathrm{i \in I_{n}}} \left( \mathrm{c_{c,i}} S_\mathrm{{i}} + \mathrm{c_{f,i}} C_{c,i} + \sum_{\mathrm{t \in T}} \mathrm{c_{v,i}} \sum_\mathrm{r \in R}  X_\mathrm{{out,i,r,t}} \right)
\end{equation}

Network costs $C_\mathrm{netw}$ are treated similarly. Again, the investment costs for existing networks are zero ($C_{c,l} = 0$. For new networks, $C_{c,l}$ is size and distance dependent and the definition can be found in the Supplementary Information.
\begin{equation}
    \label{eq:Cost_Networks}
    C_\mathrm{netw} = \sum_\mathrm{{g \in G}} \sum_{\mathrm{l \in L_{g}}} \left( C_{c,l}  + \mathrm{c_{f,l}} C_\mathrm{c,l} + \sum_{\mathrm{t \in T}} \mathrm{c_{v,l}} F_\mathrm{l,t} \right)
\end{equation}

Import costs encompass costs for electricity and natural gas imports from beyond the system boundaries:
\begin{equation}
    \label{eq:Cost_Imports}
    C_\mathrm{imp} = \sum_{\mathrm{n \in N}} \sum_{\mathrm{r \in R}} \sum_{\mathrm{t \in T}} \left( \mathrm{c_{imp,r}}  X_\mathrm{{imp,n,r,t}} \right)\\
\end{equation}

Lastly, carbon costs account for the cost of carbon emissions of imports and conversion technologies:
\begin{equation}
    \label{eq:Cost_Emissions}
    C_\mathrm{CO2} = \mathrm{c_{CO2}} (E_\mathrm{tec} + E_\mathrm{imp})
\end{equation}

\textbf{Cost minimization at emission reduction target}. In this case, the objective function is the same as for the cost minimization problem (eq. \ref{eq:MinCost}), but a constraint on total annual emissions is formulated:
    \begin{equation}
        \label{eq:MinCostAtTarget}
        E_\mathrm{tec} + E_\mathrm{imp} \leq \overline{\mathrm{E}}
    \end{equation}

\subsubsection{Constraints}
Constrains fall into three different types: (i) energy balance per node $n$, carrier $r$ and time-step $t$, (ii) technology performances and (iii) network performances. Eq. \ref{eq:EnergyBalance} depicts the energy balance accounting for demand, technology inputs and outputs, network inflows and outflows as well as imports and exports beyond the system boundaries.

\begin{equation} \label{eq:EnergyBalance}
    \begin{split}
        \mathrm{D_{n,t,r}} = & \sum_{\mathrm{i \in I_{n}}} X_\mathrm{{out,i,r,t}} - X_\mathrm{{in,i,r,t}} \\
        & + \sum_\mathrm{{g \ in G}} \sum_{\mathrm{l \in L_{g,n}}} F_\mathrm{{out,l,r,t}} - F_\mathrm{{in,l,r,t}} - F_\mathrm{{cons,l,r,t}}\\
        & + X_\mathrm{{imp,n,r,t}}\\
        & \forall n \in N, r \in R, t \in T\\
    \end{split}
\end{equation}

The imports of electricity are further limited at each node to the interconnection to neighboring countries that are not modeled (eq. \ref{eq:ImportConstraint}). Exports are not possible. The natural gas pipeline network is not part of the optimization model, and thus imports of natural gas to each node are unlimited.

\begin{equation}
    \label{eq:ImportConstraint}
    X_\mathrm{{imp,el,r,t}} \leq \overline{\mathrm{X}}_\mathrm{{imp,el,r}} \ \forall t \in T
\end{equation}

The technology and network models can be found in the supplementary information. The model is formulated formulated with the python based energy system modeling software AdOpT-NET0 \cite{Wiegner2024}. It relies heavily on the modeling language pyomo and the model is solved with Gurobi 10.

\subsection{Input Data}
All required input data was gathered and pre-processed to resemble the energy system of the modeled states in 2030 as projected by the National Trend scenario of the \ac{TYNDP} \cite{ENTSO-E2022a} and the \ac{ERAA} \cite{ENTSO-E2022}. Electricity demand, onshore renewable generation profiles and onshore installed capacities are only available on a national level and thus need to be spatially allocated to the node definition of this work. Offshore wind generation was modeled on a farm-by-farm level considering different turbine types and hub heights. A brief description of the data sources and pre-processing steps taken can be found in table \ref{tab:DataSources}. Additionally we added a detailed discussion to the Supplementary Information. The whole dataset is available for download as a Supplementary Information.

\input{Tables/00_DataSources}

\FloatBarrier
\subsection{Scenario Definition\label{sec:ScenarioDefinition}}
Table \ref{tab:ScenarioDefinition} shows an overview of all scenarios defined in this work. 

\input{Tables/00_ScenarioDefinitions.tex}

In the \textbf{Reference} scenario for 2030, we only take into account the electricity system as projected by \ac{ENTSO-E} in their National Trend scenario. Therefore, investment into additional infrastructure is not possible and only the operation of the system is cost-optimized. For 2040, electricity demand and hydrogen demand are increased, but the same starting system as for 2030 is used. However, renewable capacities (onshore wind, offshore wind and onshore PV) can be expanded and are thus optimized as well. The 2040 Reference scenario thus shows the required expansion of renewable capacities to meet additional electricity demand. In all other scenarios, we study the role of supporting infrastructure to lower annual cost and/or annual emissions (brownfield approach). In line with Wiegner et al. (2024) we study the role of three different integration measures: Electricity Transmission \textit{T-Scenarios}, Electricity Storage in Batteries \textit{S-Scenarios}, and the role of Hydrogen both as an electricity storage and for industrial applications \textit{H-Scenarios} \cite{Wiegner2024a}. For all 2040 scenarios, also renewable capacities can be expanded. 

\textbf{The role of electricity grids.} In an ideal world, electricity grids can be expanded between any two locations if it is required and cost efficient. This is the case for the first scenario \textit{T-All}. An expansion of both AC and DC lines are possible, depending on the location (generally offshore: DC, onshore: AC). However, various real-world challenges impede the unconstrained expansion of grids. These barriers include public resistance to new onshore lines, technical or market barriers for implementing DC meshed grids offshore, or political constraints on cross-border connections. To study the effects of these challenges, we defined three additional scenarios: Only onshore expansions (\textit{T-1}), only offshore expansions (\textit{T-2}) and no additional transmission lines across national borders (\textit{T-3)}.

\textbf{The role of electricity storage in batteries.} Electricity storage is seen as an important measure to match supply and demand inter-temporally. We study the overall effect of electricity storage in the scenario (\textit{S-All}), where we allow Li-Ion Batteries to be installed at all onshore and offshore nodes. However, we limit the installed capacity at offshore nodes to 140 GWH per offshore node. This is approximately equal to the weight capacity of two large offshore platforms. Additionally, we study the role of offshore and onshore storage individually as both developments might be hindered by regulation, market design or public/political opposition. \textit{S-1} allows for onshore storage only, \textit{S-2} for offshore storage only. Lastly, we defined an additional scenario \textit{S-All-HP}, with a higher power-to-energy ratio of the energy storage. This is to identify if discharging power or energy capacity is the driving factor behind the other results.

\textbf{The role of hydrogen.} Hydrogen is set to play a vital role in the energy transition - in the electricity sector for energy storage as well as for other applications (industrial, transport, heating, etc.). We study the overall role of hydrogen in the \textit{H-All} scenario, in which blue hydrogen can be replaced by hydrogen from carbon-free sources (nuclear or renewable). Additionally, hydrogen can serve as a storage medium that can be reconverted to electricity via existing gas turbines or new fuel cells. The re-conversion in existing gas turbines is possible with ad-mixing hydrogen to natural gas up to 5\% of the total energy input. Additionally, we investigate the trade-offs between onshore and offshore hydrogen production by restricting electrolysis to onshore locations only (\textit{H-1}) and to offshore locations only (\textit{H-2}). To study the importance of hydrogen storage we exclude it in scenario \textit{H-3}. Lastly, we restrict long-distance hydrogen transport by not allowing for any hydrogen networks between different nodes (\textit{H-4}). As such, hydrogen produced at one node also needs to be stored or consumed at the same node.

\textbf{Synergies of all three integration measures.} In addition to the individual role of the three integration measures mentioned before, we define one scenario, that looks at the synergies of all measures. In this scenario (Scenario \textit{Synergies}), all technologies  and network expansions are allowed, serving as a global optimum or reference case.

\textbf{Common Assumptions.} Additionally, all scenarios have the following assumptions in common:
\begin{itemize}
    \item Electricity import from outside the system boundaries is possible, but discouraged by imposing a high import price (1000 EUR/MWh) and a high emission factor (0.8 t/MWh$_el$ equal to an inefficient coal fired power plant).
    \item The annual hydrogen demand in 2030 is assumed to be 370TWh (and 757TWh in 2040) is equally distributed over all timesteps. The allocation method of hydrogen demand to each node can be found in the Supplementary Information. 
    \item If not supplied by electrolysis from within the system, hydrogen is produced by steam methane reforming with carbon capture at the cost of natural gas (40 EUR/t + emission costs) and its respective emission factor (0.108 t/MWh).
    \item The assumed carbon price is 80 EUR/t in 2030 and 100 EUR/t in 2040, applicable to the combustion of natural gas, coal, and oil as well as electricity imports.
\end{itemize}

\FloatBarrier\newpage
\section{Results}
\subsection{Reference Scenario (2030)}
In the \textit{Reference} scenario, around 82 \% of electricity demand is met with renewable generations. The optimization results are in line with national goals for renewable generation of the respective countries (see Table \ref{tab:Results_Reference}). Electricity imports from outside the system play a very limited role as they only account for 0.008 \% of total electricity supply and 0.8\% of total costs. Overall, around 12.1\% of renewable generation is curtailed. This is slightly less then electricity generation from fossil fuels (13.2\%) suggesting a high emission reduction potential of flexibility measures. Norway, Denmark and the UK have a very high share of renewable generation (100\%, 97.3\% and 84.2\% respectively). In Norway this is mainly due to electricity generation from its vast hydro electric resources. In Denmark and the UK the available generation from renewable resources (mainly onshore and offshore wind) are larger than the annual sum of electricity demand. As a result, Denmark and the UK can also export significant shares of their domestic generation to other countries ($\geq$ 10\%). Belgium, on the other hand, has the largest deficit with electricity imports accounting for around 10\% of its total supply, mostly from the UK and the Netherlands. Germany is also a net importer, with most of its imports coming from Denmark during hours with excess renewable generation flowing towards the two southern nodes (DE1 and DE4). A similar pattern occurs in the UK, where the south-western node (UK4) has only limited renewable capacities and thus a large fraction of its demand is supplied by offshore wind from the north east of the UK.

\input{Tables/01_ResultsReference.tex}

The discussion of all other scenarios evolves around two main optimization objectives: First, we discuss the minimum emission point of a scenario. While this is a hypothetical case, as it would be prohibitively expensive to reduce the last possible ton of \ce{CO2}, this case estimates the reduction potential of a integration measure. Second, the cost optimal point of a scenario is discussed. If it coincides with the reference costs, it can be concluded that there is no cost reduction potential of a certain measure. Contrarily, if system costs can be reduced, the additional costs of new investments are outweighed by savings occurring in the existing system.

Subsequently, we discuss the results of the scenarios also in comparison to the \textit{Synergies} scenario. The \textit{Synergies} scenario thereby shows the maximal possible emission or cost reduction potential if all supporting infrastructure can be build. The additionally installed capacities are reported in the supplementary information and the resulting energy system design can also be visualized on the provided visualization web app. 

\FloatBarrier
\subsection{The role of electricity grids towards 2030}

\begin{figure}[h]
    \centering
    \includegraphics[width=0.98\textwidth]{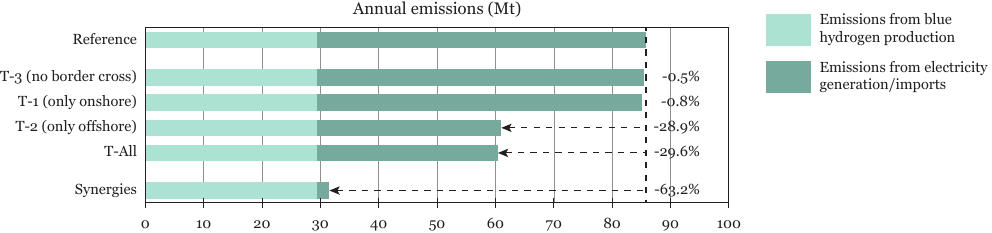}
    \caption{Emission reduction potential of electricity grids. \textit{Reference} refers to a scenario with no grid expansions, \textit{Synergies} refers to a scenario with possible expansion of storage and grid capacities as well as hydrogen conversion, storage and transport technologies.}	
    \label{fig:2030_em_grid}
\end{figure}

Figure \ref{fig:2030_em_grid} shows the \textbf{emission reduction potential} of allowing for new transmission lines. If all corridors can be expanded (\textit{T-All}) the emission reduction potential is about 30\% of the reference emissions. Removing the options to build transmission lines offshore, decreases the emission reduction potential from 30\% to less then 1\%. As all offshore lines also cross national borders, the result for not allowing for border crossings also reduces the emission reduction potential to less than 1\%. Therefore, offshore cross-border transmission lines are essential to fully realize the emission reduction potential of electricity transmission.

\begin{figure}[h]
    \centering
    \includegraphics[width=0.98\textwidth]{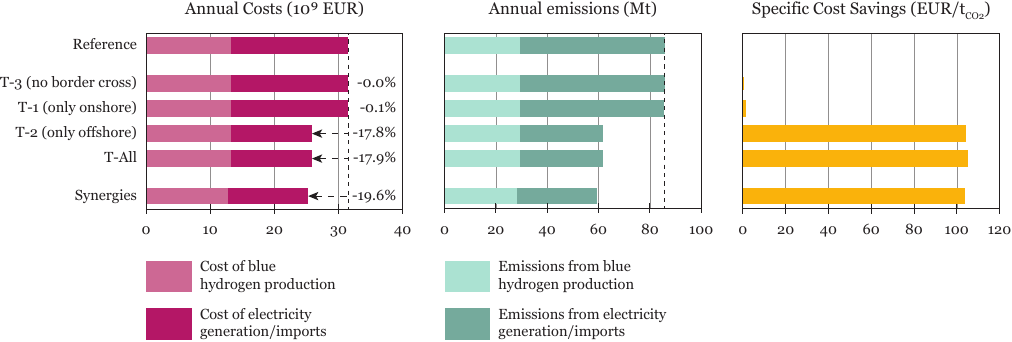}
    \caption{Cost reduction potential for the expansion of electricity grids. The figure also shows resulting emission reductions and specific cost savings per ton of \ce{CO2}. \textit{Reference} refers to a scenario with no grid expansions, \textit{Synergies} refers to a scenario with possible expansion of storage and grid capacities as well as hydrogen conversion and storage technologies.}	
    \label{fig:2030_cost_grid}
\end{figure}

Figure \ref{fig:2030_cost_grid} shows the results for the \textbf{economic potential} for all transmission scenarios. Similar to the emission reduction potential, large cost reductions are only possible if expansion corridors across borders and across the North Sea are available (\textit{T-1 (only offshore) and T-All}). The specific cost savings of these two scenarios are around 105 EUR/t. In these two scenarios, new connections are mainly made offshore between Norway and the other countries, suggesting that Norway's large hydro resources can enable large scale emission and cost reductions for other countries. Some of these new international offshore connections are also build either between offshore wind farms or between offshore wind farms and onshore nodes. As a result, some of the park-to-shore cables also serve as part of an interconnector between two countries and results in some of the wind farms being able to serve two onshore nodes. Notably, the cost-optimal solutions are very close to the emission-optimal solutions, suggesting that expanding transmission offshore can significantly reduce system costs \textit{and} emissions. Cross-border onshore connections cannot contributes much to additional cost reductions (0.1\%). This additional small reduction is possible by expanding transmission corridors between the coast and load centers further inland. 

In case only onshore transmission is allowed (\textit{T-1 (only onshore)} or the grid expansions are restricted to national projects (\textit{T-3 (no border cross)}), the cost reduction potential drops to close to zero.

Thus, the main findings in this section are:
\begin{enumerate}
    \item Investing in interconnections of countries across the North Sea is a no-regret option. It can lower overall system costs as well as emissions.
    \item The electricity transmission grid as planned onshore is sufficient to integrate the planned renewable generation in 2030. However, small expansions across country borders from renewable generation capacities offshore to load centers inland can reduce overall system costs and emissions slightly.
    \item The interconnection of existing wind farms to form simple meshed grids offshore is essential for cost and emission reduction. This results in a single wind farm serving multiple onshore nodes as well as park-to-shore cables serving as part of an international interconnector. 
\end{enumerate}

\subsection{The role of electricity storage towards 2030}
\begin{figure}[h]
    \centering
    \includegraphics[width=0.98\textwidth]{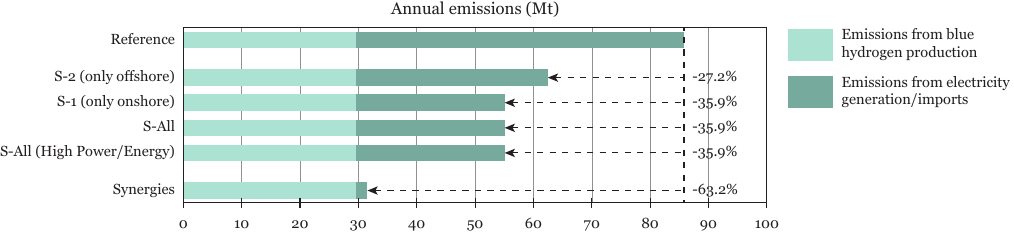}
    \caption{Emission reduction potential of electricity storage capacities. \textit{Reference} refers to a scenario with no grid expansions, \textit{Synergies} refers to a scenario with possible expansion of storage and grid capacities as well as hydrogen conversion and storage technologies. \textit{S-All (High Power/Energy)} refers to a scenario in which the fixed power-to-energy ratio is increased from 0.3 to 1 compared to \textit{S-All}.}	
    \label{fig:2030_em_storage}
\end{figure}

Figure \ref{fig:2030_em_storage} shows the \textbf{emission reduction potential} of electricity storage disregarding associated costs. It shows, that electricity storage has a significantly higher emission reduction potential then expanding electricity grids. The reason for the high emission reduction is that electricity generation from natural gas power plants can be replaced with stored electricity from renewable resources. As a result curtailed electricity in the minimum emission case in \textit{S-All} can be significantly reduced compared to the reference from 117.5 TWh (12.4\% curtailed electricity) to 34.4 TWh (3.5\% curtailed electricity). This means, that almost all renewable electricity is used. Prohibiting offshore storage (\textit{S-2 - only onshore}) has no effect on the emission reduction potential. Hence, onshore storage is most relevant in reducing emissions. Contrarily, allowing for storage only offshore significantly reduces the emission reduction potential. This is not only due to limited availability of storage locations offshore, but also because the highest energy consumption is far from shore and electricity transmission is insufficient to realize similar reduction potentials as in the onshore case. 
Increasing the fixed power-to-energy ratio from 0.3 to 1 in \textit{S-All (High Power/Energy)} does not change the findings. This suggests that large storage capacities are more important then high and quickly-available power capacities. Again, we do not report abatement costs of the emission reduction cases, as it is prohibitively expensive to mitigate the last possible ton of emissions. 

\begin{figure}[h]
    \centering
    \includegraphics[width=0.98\textwidth]{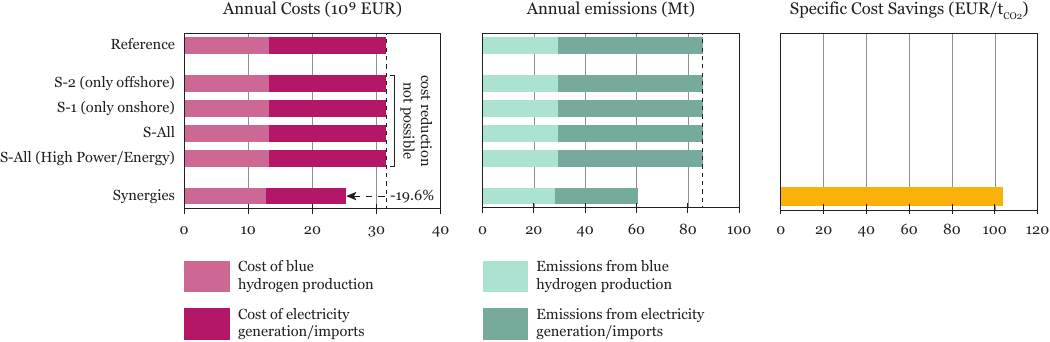}
    \caption{Cost reduction potential for the addition of electricity storage capacities. The figure also shows resulting emission reductions and specific cost savings per ton of \ce{CO2}. \textit{Reference} refers to a scenario with no grid expansions, \textit{Synergies} refers to a scenario with possible expansion of storage and grid capacities as well as hydrogen conversion and storage technologies. \textit{S-All (High Power/Energy)} refers to a scenario in which the fixed power-to-energy ratio is increased from 0.3 to 1 compared to \textit{S-All}. Note that allowing for storage additions cannot lower costs, and thus the results for all storage scenarios are the same as \textit{Reference}.}	
    \label{fig:2030_cost_storage}
\end{figure}

Allowing for electricity storage has no \textbf{cost reduction potential} (see Figure \ref{fig:2030_cost_storage}, suggesting that investment into storage is too expensive for inter-temporal balancing and it is cheaper to use flexible power plants. 

\begin{figure}[h]
    \centering
    \includegraphics[width=0.98\textwidth]{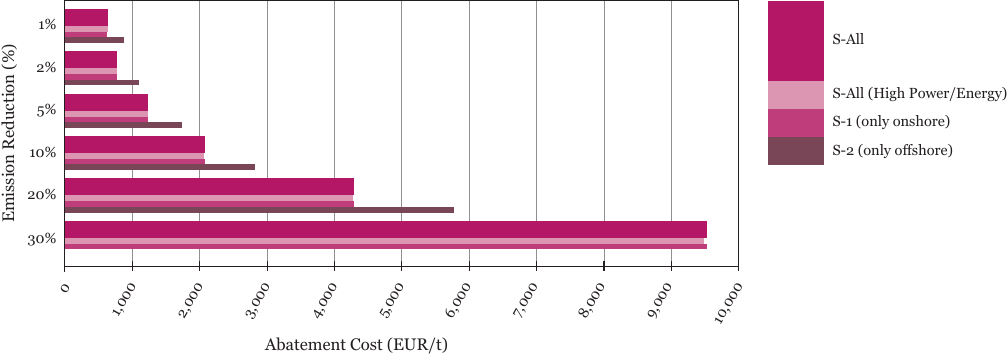}
    \caption{Abatement costs for different emission reduction targets for the addition of electricity storage capacities. \textit{S-All (High Power/Energy)} refers to a scenario in which the fixed power-to-energy ratio is increased from 0.3 to 1 compared to \textit{S-All}. Note that an emission reduction of 30\% is impossible to reach in scenario \textit{S-2 (only offshore)} and the bar is not shown respectively.}	
    \label{fig:2030_abatement_storage}
\end{figure}

To get an idea of the cost of storage to lower emissions, we report abatement costs for different emission reduction targets (see Figure \ref{fig:2030_abatement_storage}). The results suggest, that while energy storage has a high emission reduction potential, it is very costly to reach high emission reduction targets with it. For a 1\% reduction the abatement costs reach already 635 EUR/t for the scenario \textit{S-All}.  The required storage capacity to reach a 10\% emission reduction is 275 GWh for \textit{S-All}/ \textit{S-2 only onshore} and 301 GWh for \textit{S-1 only offshore}. Interestingly, the storage capacity required in the offshore only scenario is larger than in the onshore only scenario for the same emission reduction level. This means that onshore storage would still have a higher emission reduction potential, even if onshore and offshore storage costs were equal. The optimal operation of the storage is rather on a day-to-day or hourly basis. The discharged and charged energy in each week of the year is almost equal with exceptions for higher emission reduction targets. As such, seasonal storage in the projected energy system in 2030 is not required as there are sufficient back-up power plants in operation. 
Summarizing the findings of this section:
\begin{enumerate}
    \item Electricity storage has a high emission reduction potential for the power sector. 
    \item Electricity storage is an expensive emission reduction measure. However, if 'free' flexibility measures, such as demand side flexibility or vehicle-to-grid is made available to the power market, this can yield in high emission reductions without any further investments. Making available around half of the battery capacities of a projected 10 million electric cars in 2030 sums to around 300 GWh of storage capacity. If optimally used, this could yield in emission reductions of around 10\%.
\end{enumerate}

\subsection{The role of hydrogen towards 2030}
\begin{figure}[h]
    \centering
    \includegraphics[width=0.98\textwidth]{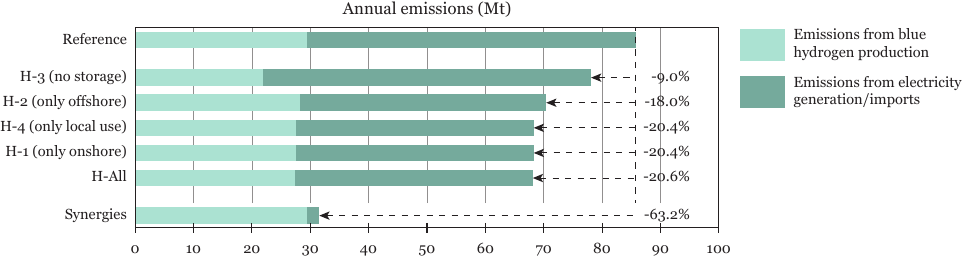}
    \caption{Emission reduction potential of hydrogen conversion, storage and transport technologies. \textit{Reference} refers to a scenario with no grid expansions, \textit{Synergies} refers to a scenario with possible expansion of storage and grid capacities as well as hydrogen conversion and storage technologies.}	
    \label{fig:2030_em_h2}
\end{figure}
The \textbf{emission reduction potential} of integrating hydrogen production, storage, transport and re-conversion into electricity into the energy system is substantial, but significantly lower than for electricity storage and electricity grid (see Figure \ref{fig:2030_em_h2}). In all scenarios allowing for hydrogen storage, hydrogen primarily serves as a storage medium and is reconverted to electricity. Similar to the battery scenarios, this integration measure can yield emission reductions by inter-temporal balancing variable renewable electricity production However, the emission reduction levels reached are lower than those achievable with batteries due to a significantly lower round-trip efficiency of hydrogen storage. The emission reduction potential for \textit{H-3 (no storage)} is only about half of the reduction potential in scenarios with hydrogen storage as carbon-neutral hydrogen can only be used to replace emission intensive blue hydrogen. This underscores the importance of hydrogen storage in maximizing the emission reduction potential compared to its direct use and as an energy carrier for transport. The two scenarios \textit{H-1 (only onshore)} and \textit{H-4 (only local use)} converge to almost the same solution as the \textit{H-All} scenario, suggesting that hydrogen transport and offshore hydrogen production offer no additional emission reduction potentials towards 2030. 

\begin{figure}[h]
    \centering
    \includegraphics[width=0.98\textwidth]{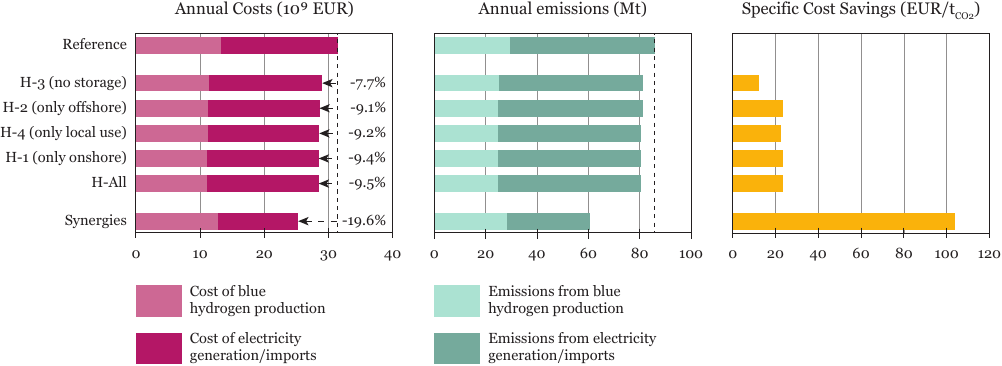}
    \caption{Cost reduction potential of hydrogen conversion, storage and transport technologies. The figure also shows resulting emission reductions and specific cost savings per ton of \ce{CO2}. \textit{Reference} refers to a scenario with no grid expansions, \textit{Synergies} refers to a scenario with possible expansion of storage and grid capacities as well as hydrogen conversion and storage technologies.}	
    \label{fig:2030_cost_h2}
\end{figure}
The \textbf{cost reduction potential} of hydrogen technologies is present, but only about half of the potential of grid expansions. However, investments in hydrogen technologies can be considered a non-regret measure that can reduce both costs and emissions. In the \textit{H-All} scenario, most hydrogen is directly utilized to replace blue hydrogen, thereby reducing both costs and emissions simultaneously. In fact, only about 4\% of total hydrogen produced is reconverted into electricity. The majority of hydrogen is produced by increasing electricity generation from nuclear power plants (57\% of all hydrogen produced), with the remainder is produced from otherwise curtailed renewable electricity. This sums to a total of 46 TWh of carbon-neutral hydrogen produced, equivalent to 16\% of projected hydrogen demand in 2030, with the rest supplied by blue hydrogen. As such, curtailment can be reduced to 6.9\% compared to 12.4\% in the reference, and the capacity factor of nuclear power plants increases from 49.1\% to 65.6\%. The required electrolysis capacities are around 15 GW with an average capacity factor of 35\%. As for the minimum emission points, the optimal solutions for  \textit{H-1 (only onshore)} and \textit{H-4 (only local use)} coincide with scenario \textit{H-All}, suggesting that long-distance hydrogen transport and offshore hydrogen infrastructure does not offer additional for cost reductions. In contrast, forcing all hydrogen to be produced offshore (scenario \textit{H-2 (only offshore)} yields lower cost reductions and reduces the emission reduction potential. This is due to higher investment costs for offshore electrolysis and a lower flexibility as hydrogen is not produced at a possible end-use location. For the same emission reduction target (e.g. 10\%) the required electrolysis capacity offshore is larger than the required onshore capacity. As such, even if the specific investment costs of offshore electrolysis was the same as onshore, reaching similar reduction targets requires higher investments.

For all scenarios, higher emission reduction targets increase the share of hydrogen stored and reconverted to electricity. While in the cost-optimal case, the share is smaller than 4\% it increases to 84\% in the emission-optimal case (scenario \textit{H-All}). For reduction targets beyond 10\%, also the installed natural gas plant capacities are not sufficient for re-conversion, and additional fuel cells are required.

We summarize the main take-aways of this section as:
\begin{enumerate}
    \item In the cost-optimal case, carbon-free hydrogen is mainly used to replace carbon-intensive blue hydrogen saving both emissions and costs.
    \item With higher emission reduction targets, the importance of hydrogen storage for intertemporal balancing grows. As such, also hydrogen storage becomes more important.
    \item Long distance hydrogen transport plays only a minor role, both for emission and cost reductions. 
\end{enumerate}

\subsection{Synergies of all integration measures towards 2030}
In the \textit{Synergies} scenario, all three integration measures are possible: new transmission and storage capacities as well as new hydrogen conversion, storage or transport assets. The results for the cost optimal and emission optimal cases are depicted in Figures \ref{fig:2030_em_grid}-\ref{fig:2030_cost_storage} and \ref{fig:2030_em_h2}-\ref{fig:2030_cost_h2}, alongside the results of the individual integration measure. These figures indicate, that no single integration measures can reach the maximum potential for emission or cost reduction with the projected renewable capacities in 2030. Therefore, synergistic benefits from multiple integration measures are evident. 

The highest \textbf{emission reduction} can be achieved through a combination of transmission line expansions and electricity storage capacities, enabling both inter-regional and inter-temporal balancing. Inter-regional balancing is preferred due to lower losses in electricity transmission compared to the charging and discharging processes of electricity storage. Hydrogen technologies, however, do not contribute to the maximum emission reduction potential. In the solution with highest emission reduction, curtailment is reduced to zero and as such all generated electricity from renewable sources is used to cover electricity demand. The required storage and transmission capacities are disproportionately high, and therefore, we do not report them here.

For the highest \textbf{cost reductions}, it is optimal to deploy a mix of transmission line expansions, electrolysis, and hydrogen storage. The electricity line expansions follow the same corridors as previously discussed, resulting in a meshed grid in the North Sea with transmission lines serving both as interconnectors and park-to-shore cables. Hydrogen serves a dual purpose: nearly all (99.3\%) of the carbon-neutral hydrogen produced is directly used to replace blue hydrogen, while the remainder is stored in hydrogen storage facilities and later mixed into the fuel for natural gas-fired power plants to generate electricity. Given the expected renewable capacities in 2030, approximately 5.4 GW$_\mathrm{el}$ of electrolysis capacity is required for the cost-optimal case, with most of it being installed in the UK and powered by a mix of nuclear and renewable electricity. As shown in Figure \ref{fig:2030_cost_grid}, the cost difference between the \textit{T-All} and \textit{Synergies} scenarios is only 1.7 percentage points. The full potential is realized when, in addition to expanding electricity grids, hydrogen technologies are also deployed. However, given the small difference, the role of hydrogen technologies in achieving the full cost reduction potential is limited, assuming that grid expansions can be successfully implemented.

\begin{figure}[h]
    \centering
    \includegraphics[width=0.98\textwidth]{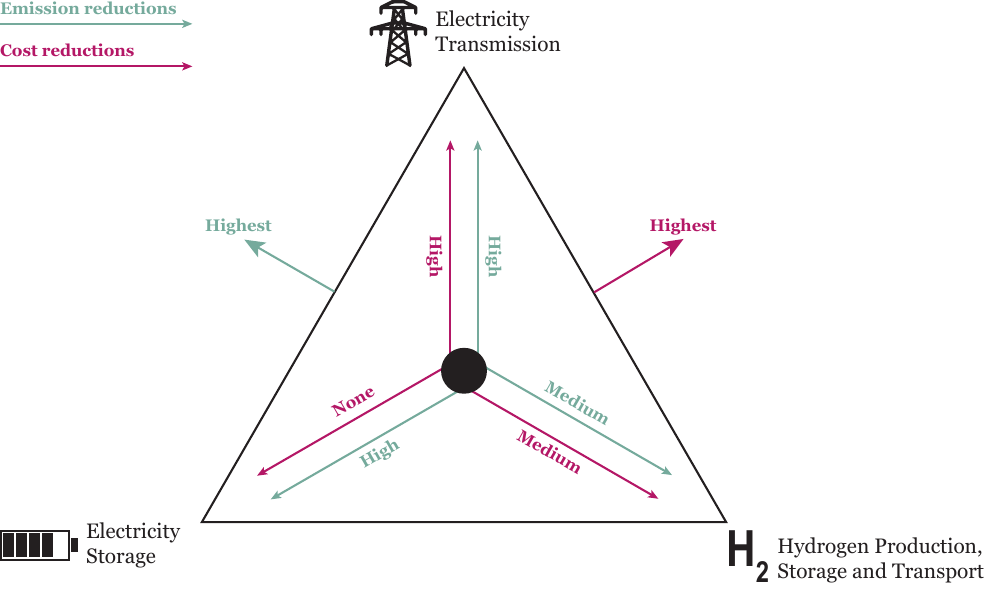}
    \caption{Qualitative summary of the findings of the role of supportive infrastructure for the North Sea region in 2030. The middle of the triangle depicts the reference scenario with no additional infrastructure being allowed to build. The corners of the triangle depict the three integration measures studied with an indication of their cost reduction and emission reduction potential.}	
    \label{fig:2030_summary}
\end{figure}

The main findings can be summarized as follows (also summarized in Figure \ref{fig:2030_summary}):
\begin{enumerate}
    \item The full potential for emission reduction can only be realized through the expansion of both transmission grids and electricity storage capacities.  Hydrogen technologies cannot contribute to reach the maximal possible emission reduction.
    \item The largest cost reductions are achieved through the expansion of the transmission grid offshore. A small, additional reduction can be gained by producing carbon-neutral hydrogen to replace blue hydrogen.
    \item Expanding the electricity grid across the North Sea emerges as a no-regret strategy for reducing both costs and emissions in the region.
\end{enumerate}

\subsection{Pathway towards 2040}
For the 2040 case, we assumed the same transmission, generation, and hydro storage capacities as in 2030. However, electricity and hydrogen demand increases according to the projections from \ac{TYNDP}. To meet this additional demand, renewable energy capacities can be expanded within feasible limits. The \textit{Reference} scenario for the year 2040 now refers to a case in which only \ac{vRES} assets can be extended (i.e. onshore wind, offshore wind and photovoltaic). For new offshore wind parks, also the respective park-to-shore cable is allowed, but no interconnection to other onshore nodes or wind farms is permitted. The other scenario definitions remain unchanged from 2030.

The \textbf{emission reduction potential} is very high and it is possible to reduce emissions of the electricity sector almost to zero (2.1 Mt) in the \textit{Reference} scenario. This can be achieved by an extreme expansion of renewable capacities. The same is the case for all other scenarios, but emission from hydrogen production remains. If all hydrogen technologies are allowed, emissions from both hydrogen production and electricity generation can reduced to zero. However, also in this case, installed capacities are unrealistically high and we thus do not discuss these cases further.

\begin{figure}[h]
    \centering
    \includegraphics[width=0.98\textwidth]{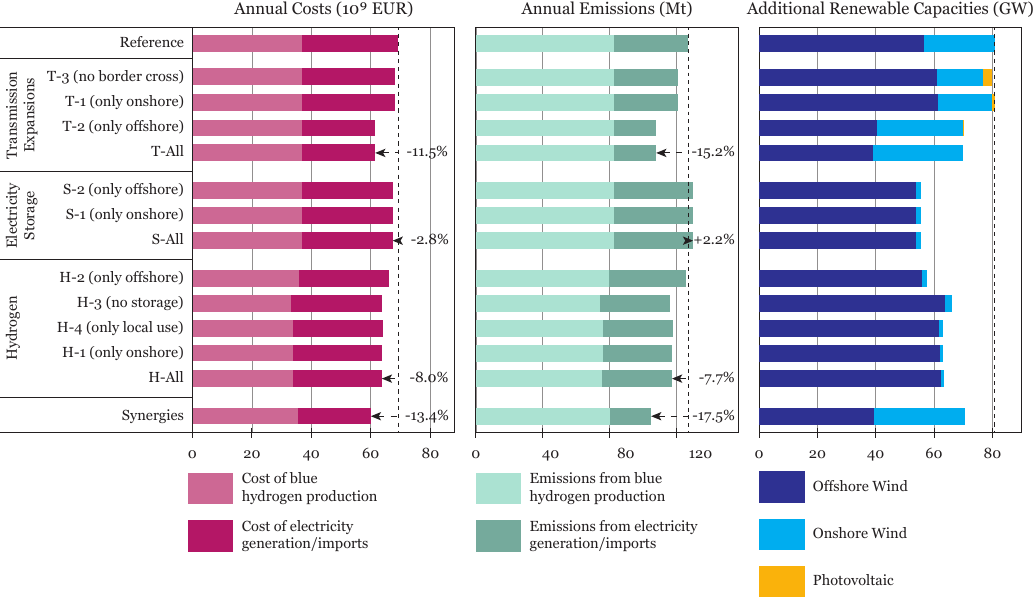}
    \caption{Cost reduction potential of hydrogen conversion, storage and transport technologies. The figure also shows resulting emission reductions and specific cost savings per ton of \ce{CO2}. \textit{Reference} refers to a scenario with no grid expansions, \textit{Synergies} refers to a scenario with possible expansion of storage and grid capacities as well as hydrogen conversion and storage technologies.}	
    \label{fig:2040_cost}
\end{figure}

The \textbf{cost reduction potential} across all scenarios is shown in Figure \ref{fig:2040_cost} alongside total emissions and new capacities for \ac{vRES} assets. The installed capacities of the supportive infrastructure are reported in the supplementary information and the resulting energy system design can also be visualized on the provided visualization web app. 

In the \textit{Reference} scenario, approximately 80 GW of additional \ac{vRES} capacity is identified as cost-optimal to meet the increased demand compared to 2030. These additional renewable capacities not only meet the new demand but also help reduce emissions in the electricity sector by around 20\% compared to 2030. Offshore wind makes up the majority (70\%) of the new \ac{vRES} installations, with additional onshore wind installations in the UK and Norway. These new capacities are distributed across the North Sea region. It is important to note that total emissions in 2040 are higher than in 2030 due to increased blue hydrogen production, making a direct comparison of total emissions between the two years inapplicable.

As with the 2030 case, electricity grid expansions across the North Sea emerge as a no-regret option (\textit{T-All}). These expansions can significantly reduce (a) the required investment in renewable assets by 18\%, (b) total system costs by 12\%, and (c) total emissions by 17\%. The expansion corridors are similar to those optimal in 2030, primarily connecting offshore wind farms with onshore nodes across borders. In scenarios where cross-border transmission lines and international connections are not allowed (\textit{T-3} and \textit{T-1}), the economic and environmental benefits of grid expansion are significantly reduced. Hence, offshore cross-border connections remain essential in 2040 for realizing economic and environmental gains.

In contrast to the 2030 case, deployment of electricity storage can reduce costs slightly. This is due to the possibility to expand renewable capacities, and the additional demand in 2040 can be covered by expanding \ac{vRES} or with storage capacities. As a result, investment in electricity storage can significantly lower required investment into \ac{vRES} capacities compared to the \textit{Reference} scenario. The required storage capacity to reach the minimal costs is about 5 GWh and the storage optimally operates as a short term storage in the range of hours to days. The solution for the scenarios \textit{S-All} and \textit{S-1 (only onshore)} is equivalent with all storage being installed onshore. Allowing for storage only to be installed offshore increases the required \ac{vRES} capacities slightly and results in slightly higher total annual costs. Notably, in the cost optimal scenario, electricity storage can only reduce system costs and required additional \ac{vRES} capacities. However, operating the storage systems cost optimally does not reduce emissions, as \ac{vRES} contribute less to electricity supply and as a consequence generation from fossil fuel power plants increase compared to the \textit{Reference}.

Unlike in 2030, electricity storage can help to reduce total system costs slightly. Given the higher electricity demand and the possibility to increase \ac{vRES} capacities, storage capacities can complement \ac{vRES} expansions. Possible investment in electricity storage thus reduces the required \ac{vRES} capacities compared to the \textit{Reference} scenario. The optimal storage capacity is about 5 GWh, and it operates as short-term storage, with cycles typically ranging from a few hours to days. The solution for the scenarios \textit{S-All} and \textit{S-1 (only onshore)} is equivalent with all storage being installed onshore. As such, offshore storage leads to slightly higher costs and increased optimal \ac{vRES} capacities. Notably, while storage reduces system costs and \ac{vRES} capacity needs, it does not lower emissions, as the use of fossil fuel plants increases when \ac{vRES} contributions decrease.

Hydrogen technologies towards 2040 have a moderate cost reduction potential, but as the other two integration measures, they can significantly reduce the required additional \ac{vRES} capacities. In the \textit{H-All} scenario, 90\% of produced hydrogen is directly used to replace blue hydrogen and thus can lower emissions from blue hydrogen production. The remaining 10\% are used for inter-temporal balancing being reconverted in existing gas power plants and newly build fuel cells. Hydrogen also serves as a transport carrier transporting hydrogen along a north-south connection in the UK to bring hydrogen from offshore wind in the north of the UK to the load centers in the South. If hydrogen storage is not allowed (\textit{H-3 (no storage)}), the pipeline network further extents towards Norway and Denmark along corridors re-purposing existing natural gas pipelines. However, in none of the scenarios, it is cost efficient to expand \ac{vRES} capacities to produce hydrogen. The cost optimal electrolysis capacities are between 10GW (\textit{H-2, only offshore}) to around 30GW (for all other hydrogen scenarios) and are distributed over Germany, Denmark, The Netherlands and the UK. As in the 2030 case, the produced hydrogen comes from a mix of renewable and nuclear electricity.

Hydrogen technologies have moderate cost reduction potential by 2040, but, similar to other integration measures, they significantly reduce the need for additional \ac{vRES} capacities. In the \textit{H-All} scenario, 90\% of the hydrogen produced replaces blue hydrogen, thereby reducing total system emissions. The remaining 10\% is used for balancing, reconverted into electricity via gas power plants and new fuel cells. Hydrogen also functions as a transport medium, with a north-south pipeline in the UK delivering hydrogen from offshore wind parks in the north to demand centers in the south. 
If hydrogen storage is not allowed (\textit{H-3}), the pipeline extends towards Norway and Denmark, repurposing existing natural gas pipelines. Optimal electrolysis capacities range from 10 GW (\textit{H-2, only offshore}) to around 30 GW in all other hydrogen scenarios, distributed across Germany, Denmark, the Netherlands, and the UK. As in 2030, hydrogen production is based on a mix of renewable and nuclear energy sources. Prohibiting onshore hydrogen production (\textit{H-2, only offshore}) significantly reduces both cost and environmental benefits.

In the \textit{Synergies} scenario, investments in electricity storage and hydrogen technologies as well as grid expansions are allowed. Total cost savings amount to 13.4\%, with grid expansions and hydrogen storage being the drivers. Hydrogen transport and electricity storage are not cost-effective in this scenario. Notably, the difference between \textit{T-All} and \textit{Synergies} is marginal (1.9 percentage points), indicating that most savings can be achieved through grid expansions, with hydrogen production contributing the remaining 1.9 percentage points. Most hydrogen (98\%) is used to replace blue hydrogen production, while the remaining 2\% is stored and reconverted into electricity in natural gas power plants. Investment in new \ac{vRES} capacities represents the largest portion (68\%) of total investments into new infrastructure.

\textbf{Comparison with 2030.}
To cover additional electricity demand, we allowed for \ac{vRES} capacity expansions in 2040. However, the expansion of \ac{vRES} capacities also opens a new pathway to reduce costs and emissions. The importance of the three integration measures, however, remains consistent with 2030: grid expansions across the North Sea offer the highest cost-reduction potential, while hydrogen technologies and electricity storage provide lower economic benefits. The cost-optimal mix includes large grid expansions and smaller electrolysis capacities, with hydrogen primarily used directly, rather than for storage or transport.

The key difference between the year 2030 and 2040 is the role of electricity storage: Co-optimizing the capacities for additional renewable capacities and electricity storage makes relatively small storage capacities cost competitive: Approximately 5.5 GWh of electricity storage reduces the required \ac{vRES} additions by around 25.6 GWh.

The main findings for the year 2040 are thus as follow:
\begin{enumerate}
    \item The three integration measures remain similarly effective in 2030 and 2040: Grid expansions provide the highest economic and environmental benefits with hydrogen production to replace blue hydrogen plays a smaller role.
    \item All three integration measures can reduce the required \ac{vRES} capacity to meet increased demand.
    \item Total system emissions (from hydrogen production and electricity generation) can theoretically be reduced to zero in 2040 through large-scale \ac{vRES} additions.
\end{enumerate}

\FloatBarrier\newpage
\section{Discussion}
\subsection{Results in light of the current governance framework}
Our work has studies the energy system of the North Sea region from a global optimization point of view. This perspective does not align with real-world decision processes and policy making. Achieving a global optimum is thus unlikely, given that central planning is absent and policy making does typically not consider a trans-boundary energy system. Instead, decisions are often made from the perspective of individual actors within a single jurisdiction (governments, administrative bodies, investors), whose priorities often diverge from what is beneficial for the system as a whole. Consequently, actual policy making is prone to diverge from a cost- or environmentally effective path due to lobbying and national interests \cite{Pezzey2014, Suehlsen2014}. Additionally, highly relevant stakeholders, such as Transmission System Operators (\acs{TSO}), have limited lobbying power - given their government oversight and lack of direct financial incentives - and others with greater influence may push agendas that do not align with a global energy system perspective \cite{Biancardi2021a}. As such, current governance system is not well-equipped to make decisions that are best for the system as a whole, particularly in the context of trans-boundary and trans-sectoral collaboration.

While investment into electricity storage, hydrogen production, storage and use can be based on rather small, project-based planning processes. Grid expansions across the North Sea are fundamentally different and highly complex. This complexity arises because grid expansions require large-scale cooperation across different actors, including \acs{TSO}s, spatial planning agencies, national governments, EU bodies, and offshore wind park operators, not to mention the different legal systems between EU and non-EU countries \cite{Platjouw2018, Jay2016, QueroGarcia2019, Kusters2023}. Moreover, such transnational projects often have asymmetric benefits, leading to potential conflicts between market parties and between countries \cite{GorensteinDedecca2018}. To overcome these challenges, harmonized rules and cross-boundary governance bodies are necessary for an effective planning and implementation process. Furthermore, to ensure equitable distribution of benefits and to mitigate public opposition effective compensation payments between actors need to be agreed upon. 

While this study has focused primarily on the affordability and environmental effects of integration measures, it is important to acknowledge that other considerations are also critical in decision-making. Energy security, employment, geopolitical interests, access to raw materials, social acceptance, and life-cycle impacts of technologies all play a significant role. For instance, while grid expansion may offer significant system benefits, it may also raise concerns related to energy security or geopolitical dynamics, especially given the complex inter-dependencies it creates between countries. Moreover, the social acceptance of large infrastructure projects cannot be overlooked. Even if such projects are economically and environmentally beneficial, they may face opposition due to perceived negative impacts on local communities or the environment. This opposition can stem from concerns about landscape changes, potential harm to local ecosystems, or disruptions to existing industries. Addressing these concerns through inclusive planning processes and ensuring that benefits are equitably distributed are essential for gaining public support.

Another important aspect is the distinction between operational expenditures of existing infrastructure and investment into new infrastructure. Investing in new infrastructure, is capital-intensive and raises the question of who will bear these costs. While a global system perspective suggests that such investments are necessary for achieving long-term sustainability and cost efficiency, the upfront financial burden is substantial and often cannot be financed out-of-pocket. To bridge this gap, innovative financing mechanisms and policy frameworks are required to distribute costs fairly among stakeholders, including governments, private investors, and end-users. Additionally, there may be a need for public subsidies or incentives to offset the high initial costs and to encourage investment in projects that provide long-term benefits to the entire system.

\subsection{Model Limitations}
Firstly, we did not run the scenarios for \textbf{different climate years}, and as such our results are based on the weather patterns of one typical year only. Although the figures would change for different climate years, the variability of \ac{vRES} remain and thus we expect the general trends and conclusions to also hold for different climatic years.

\textbf{Technology costs} were fixed and assumed for the year 2030. Different cost assumptions, especially of less mature technologies such as meshed \ac{HVDC} grids would change the cost reduction potential respectively, while the emission reduction potential of the 2030 scenarios remain. However, we do not anticipate that changes in costs will reverse the trend observed in this study. As the difference in the scenarios are large, the same trends observed would also hold for different cost assumptions

The \textbf{spatial aggregation} in our model does not capture grid congestion at the regional level or of the distribution grid. This limitation means that local grid congestion might be overlooked. However, the focus here is on larger-scale transmission grid expansions, where significant system-wide benefits can be realized. Investment in distribution grids is generally recognized as necessary, regardless of our findings.

We assume average values for the \textbf{costs and emission factors of blue hydrogen}. If blue hydrogen was available at a different price, the economic potential of carbon-free hydrogen as a replacement would change respectively. However, the emission reduction potential would remain unchanged, as all hydrogen in our model is used for energy storage in the minimum emission case. The role of hydrogen storage for emission reduction, however, might decrease if blue hydrogen had a higher emission factor (e.g. being produced from natural gas without \ac{CCS})

\textbf{Technological maturity}. We assume that all integration measures considered in the model can be readily implemented. In reality, large-scale electrolysis, hydrogen storage and High Voltage Direct Current (HVDC) transmission in a meshed grid are not yet technologically mature or widely deployed. Although these technologies are considered feasible and crucial for future energy systems, their current lack of maturity means that our results should be interpreted with caution, particularly concerning their near-term applicability.

Finally, our study does not consider the role of \textbf{energy storage for short-term balancing}. While we focus on the economic and emission reduction potential of energy storage, it is important to acknowledge that energy storage can also play a critical role in maintaining system stability. 

\FloatBarrier\newpage
\section{Conclusion}
This study has studied the role of supportive infrastructure in integrating high shares of variable renewable energy sources (\ac{vRES}) in the North Sea region, focusing on the year 2030 and 2040. For 2030, we have taken the expected capacities of electricity generation and transmission as well as electricity and hydrogen demand as exogenous and allowed for new assets of the following three integration measures: (i) and expansion of the electricity grid, (ii) investment into electricity storage and (iii) investment into electrolysis, hydrogen transport and storage. For 2040, we considered the same \ac{vRES} as in 2030 but with the possibility to expand these. We assessed both the emission reduction potential and the economic benefits from a system cost perspective for each integration measure. Below, we address the research questions posed in this work.

\textit{How can supportive infrastructure (grid expansions, electricity storage and hydrogen technologies) contribute to welfare gains and emission reduction in the short (2030) and medium term (2040)?}

Grid expansion across the North Sea, interconnecting countries and wind farms, offers the highest potential for both emission and cost reductions in both 2030 and 2040. Hydrogen technologies can complement the expansion of the grid by offering an additional small economic benefit. Electricity storage is not cost-effective when competing with hydrogen production and grid expansion. However, large-scale deployment of electricity storage can lead to significant emission reductions, albeit at high costs.

\textit{Are there synergies between the three integration measures and does e.g. one enable another becoming beneficial?}
The combination of hydrogen production and storage with grid expansions offers the highest cost reduction potential in both 2030 and 2040. While electricity storage, when combined with grid expansion, can lead to the greatest emission reductions in the short term, the storage capacities required are very large. Making existing (and free) storage capacities, e.g. demand side management or vehicle-to-grid, available for balancing is thus a no-regret option.

\textit{What is the welfare effect and the impact on the emission reduction potential if a measure is infeasible due to political, legal, social or technical constraints?}
Transmission infrastructure across the North Sea, particularly connections between Norway and other countries, is crucial to realizing the full potential of \ac{vRES} in terms of both emission and cost reductions for 2030 and 2040. While planning and building an interconnected electricity grid in the North Sea is the most challenging integration measure, due to technical, legal and governance challenges, it is also the most rewarding. Therefore, prioritizing the development of these corridors and agreeing on compensation mechanisms or fair cost sharing should be a key policy focus.

If hydrogen can be produced by electrolysis, it is most economic to directly use it as a fuel or feedstock at the place of production. Hydrogen storage, along with reconversion into electricity, can contribute to modest cost and emission reductions. In contrast, hydrogen transport is not essential in the short and medium term. The versatility of hydrogen, which can be used both for industrial purposes and power generation, makes it a valuable, though secondary, integration measure.

Overall, this study underscores the critical importance of expanding electricity grid connections across the North Sea to fully leverage the potential of \ac{vRES} in the region. While the production of carbon-neutral hydrogen from \ac{vRES} or nuclear power can also contribute to emission and cost reductions, its impact is less significant compared to grid expansion. Electricity storage, though a powerful tool for reducing emissions, is an extremely expensive option, especially in the short term. Towards 2040, electricity storage can mostly reduce the required additional \ac{vRES} capacities to meet growing electricity demand.

\newpage
\bibliographystyle{unsrt}
\bibliography{library}

\FloatBarrier\newpage
\part*{Supplementary Information}
\setcounter{section}{0} \renewcommand{\thesection}{S\arabic{section}}

\section{Input Data}
The following electricity conversion technologies are considered in this study: natural gas, coal, nuclear, hydro run-of-river, pumped hydro (reservoir), pumped hydro (closed loop), pumped hydro (open loop), biomass, solar, wind onshore and wind offshore. The national capacities are allocated to our node definition using allocation keys. The allocation methods for each technology are described in the following sections. Table \ref{tab:02_InstalledCapacities_DataSources} shows the data sources for each of the technologies.

\input{Tables/02_InstalledCapacities_DataSources.tex}

\subsection{Installed Capacities}

\subsubsection{Coal, Gas, Nuclear, Oil, Hydro}
Installed capacities on a national level for coal, gas, nuclear, oil, hydro storage and hydro run of river plants are retrieved from the \acs{TYNDP} 2022 for the scenario National Trend \cite{ENTSO-E2022a}. The website of TYNDP states: 'Other non RES include mainly CHP that is used in district heating \& industry. Fuel use can be gas, coal, lignite, and oil. The \ce{CO2} content of ONR technologies depending on the technology and have been considered into the \ce{CO2} budget.' As such, we have added its capacity to the category 'Gas'. We also added the Block 2 and 3 of the Belgian Tihange Nuclear Power Station, as it was missing in the TYNDP data. Coal powered plants in the UK were not considered, even though they were included in the TYNDP data, as coal is supposed to be phased out until the end of 2024.
The TYNDP does not provide data on separate Hydro Storage technologies, and thus national capacities were instead taken from ERAA 2022 (also from ENTSO-E) for the National Estimates \cite{ENTSO-E2022}. 

\input{Tables/02_InstalledCapacities_NonOffshore_Country.tex}

The national capacities are then allocated to each node using the spatially resolved power plant data base by Gotzens et al. (2019) \cite{Gotzens2019}. We therefore calculated the installed capacity for each node per our definition in this paper and used the capacity per node as a key to allocate the TYNDP national capacities. For hydro storage, we assume that turbine power, pump power and reservoir capacities are distributed equally as Gotzens et al. (2019) only give turbine capacities. Furthermore, Gotzens et al. (2019) do not distinguish between open loop and reservoir plants and as such we used the allocation keys from reservoirs for both the open loop and the reservoir capacities from the ERAA data.

\subsubsection{Solar and Wind onshore}
For each country, we calculate the the installed capacity for onshore wind and solar PV on a NUTS 2 level. Therefore, we collected installed capacities for 2023 from different national sources, and potentials for solar PV and onshore wind from the ENSPRESO database \cite{JointResearchCentre2021}. As the NUTS definition has changed for the UK, UKM8 and UKM9 are redistributed based on their area. Expected national capacities for 2030 were taken from \acs{TYNDP} 2022 for the scenario National Trend \cite{ENTSO-E2022a}. 
The expected national capacities were distributed among all NUTS 2 regions using the following formula. This was done for each country. $C$ refers to capacity, $P$ to potential, $R$ for each NUTS 2 region and $N$ for national.

\begin{equation}
    C_{2030, R} = \frac{P_{R} - C_{2023, R}}{P_{N}} * (C_{2030, N} - C_{2023, N}) + C_{2023, R}
\end{equation}

The installed capacities per NUTS2 region were then used to calculate an electricity production from wind and solar for each node. Below is an overview of national installed capacities, the data per NUTS region is available upon request.

\subsubsection{Biomass}
Similar to the previous sections, we allocated the \acs{TYNDP} 2022 data for the scenario National Trend \cite{ENTSO-E2022a} to the nodes using different national sources. See the table below for a data sources and notes on allocations.

\input{Tables/02_InstalledCapacities_BiomassAllocation.tex}

\subsubsection{Wind offshore}
Capacities for all wind farms in the North Sea and the Baltic Sea were collected from several national sources that are listed below that were compared to \cite{4COffshore2023,NorthSeaEnergy2023,EuropeanCommission2023}. Only wind farms that are expected to be operational by 2030 were taken into account for this study. Additionally, wind farms that are located in the Baltic Sea or that are not expected to be connected to any other wind farm or node were added to the production profile onshore.

\input{Tables/02_InstalledCapacities_OffshoreWindSources.tex}

\subsection{Demand}
\subsubsection{Electricity Demand}
Hourly electricity demand profiles for each node are derived from national electricity demand profiles from the \acs{TYNDP} National Trend scenario for the climate year of 2008 \cite{ENTSO-E2022a}. We follow the following steps to derive a higher spatial resolution:

\begin{enumerate}
    \item For all demand time series, values below 0.1 were identified as faulty and interpolated from the previous and subsequent time step accordingly.
    \item For countries with multiple bidding zones (Denmark and Norway), demand time series were aggregated to the national level.
    \item We calculate the annual electricity demand per node $D_{tot, node}$ as a fraction of total annual national electricity demand. We therefore use allocation keys derived from PyPSA Europe published by Neumann et al. (2023) \cite{Neumann2023}.
    \item We assume that the variation in demand over time depends on the industrialization of a region: electricity demand from industry is assumed to be less variant over time then residential demand. We therefore split total annual demand into (1) industrial demand $D_{ind, national}$ and (2) other demand $D_{other, national}$ using the share of industrial electricity demand in 2021. We therefore use electricity demand data from EUROSTAT \cite{Eurostat2023a,Eurostat2023}:
    \begin{equation}
        D_{ind, national} = s_{ind} D_{tot, national}
    \end{equation}
    \begin{equation}
        D_{other, national} = D_{tot, national} - D_{ind, national}
    \end{equation}
    Note that we implicitly assume  that the share of industrial electricity demand will not change until 2030.
    \item Other, non-industrial national demand is calculated by subtracting the industrial demand profile from the original demand profile.
    \item The national industrial demand acquired are allocated to each node using keys calculated from employment statistics published by Eurostat for 2018 \cite{Eurostat2023b}. The allocation key for each node is the fraction of employees in the manufacturing sector at each node. These keys are calculated individually for each country:
    \begin{equation}
        D_{ind, node} = \frac{N_{manufacturing, node}}{N_{manufacturing, national}} * D_{ind, national}
    \end{equation}
    \item Non-industrial demand per node is then calculate from the total demand at each node:
    \begin{equation}
        D_{other, node} = D_{tot, node} - D_{ind, node}
    \end{equation}
    \item We generate hourly time series for each node by (1) evenly distributing annual industrial electricity demand over the year (flat demand) and (2) allocating the national non-industrial demand profile to each node according to its respective share of national non-industrial demand.
\end{enumerate}

\input{Tables/02_Demand_Country.tex}

\subsubsection{Hydrogen demand}
Annual hydrogen demand on a NUTS2 level is taken from Gross et al (2022) for their baseline scenario and allocated to each node \cite{Gross2022}. The nodal hydrogen demand was then distributed equally over all hours of a year to form a flat demand profile.

\subsection{Renewable generation profiles}
Table \ref{tab:02_DemandSupplyCountry} provides an overview over demand and supply of non-dispatchable generation. How generation profiles from non-dispatchable sources and inflows to open pumped hydro power plants are calculated is described hereafter.

\input{Tables/02_DemandSupply_country.tex}

\subsubsection{Onshore Wind and PV}
To calculate hourly generation profiles of onshore wind and solar PV for each node, we use the installed capacities at each NUTS2 region (as described before). For generation from wind turbines we used the wind speed at 100m and for solar PV the irradiance retrieved from the ERA5 database at the centroid of each NUTS2 region \cite{Munoz-Sabater2021}. A height correction is applied to recalculate the wind speeds to 110m using the following formula:

\begin{equation}
    \label{eq:heightCorrection}
    ws_{110m} = ws_{10m} \times \left(\frac{110}{100}\right)^{1/7},
\end{equation}

The capacity factor for each NUTS 2 region is computed based on the power curve for a 1.5 MW turbine shown in Figure abc. This capacity factor is then multiplied by the installed capacity for the respective NUTS 2 region. The resulting time series are aggregated to determine wind generation per node. 

In the case of Norway, the installed capacity for each bidding zone retrieved from TYNDP data \cite{ENTSO-E2022a} is multiplied by the corresponding capacity factor extracted from the ERAA database \cite{ENTSO-E2022}. The generation profiles of individual bidding zones are subsequently aggregated to formulate a comprehensive national generation profile.

\subsubsection{Offshore Wind}
To calculate hourly generation profiles for offshore wind farms, we used the installed capacity and the power curve for the respective turbine type installed into account. Additionally, we assumed different hubheights according to the year installed, i.e. 80m for commissioned farms before 2010, 100m for farms commissioned between 2010 and 2020 and 120m for all farms installed after 2020. As for the onshore wind generation, we used wind speeds for a height of 100m retrieved from the ERA5 database at the centroid of each farm and recalculated it to the respective height using the Equation \ref{eq:heightCorrection}  with an exponent of 0.11. 

\subsubsection{Biomass}
Electricity generation from biofuels is assumed to be a flat profile with a capacity factor of 0.53 based on data for Europe \cite{Bolson2022}.

\subsubsection{Run of River}
Electricity generation profiles from run-of-river is based on the daily flows for run of river provided by the ERAA \cite{ENTSO-E2022} The daily flows were divided equally over the hours of the day and assigned to each node based on the installed capacities at the respective node.

\subsubsection{Inflow to Open Loop Pumped Storage and Reservoir}
Similar to the run of river generation, inflows into the upper reservoir for reservoir type and open loop type pumped hydro plants was calculated based on data by the ERAA \cite{ENTSO-E2022}. The weekly values were equally distributed over the hours in a week.

\subsection{Electricity Network Topology}
To define capacities of ac and dc grids between the nodes, we used a variety of sources. As a starting point, we aggregated the grid data given by Neumann et al. (2023) \cite{Neumann2023} to the node definition used in this work. For some lines, the data did not coincide well with data provided by the ERAA 2022 report and HoogespaningsNet.com \cite{Hoogspanningsnet2024}. Where this was the case, we changed it accordingly. The capacities for transmission lines between offshore wind farms and the onshore nodes are based on the same sources as the capacities for offshore wind farms (see Table \ref{tab:offshore_wind_sources}).
Interactions with countries that are out of scope is modeled with a maximum import constraint to the respective node.

\subsection{Cost Assumptions}
\textbf{Technology Costs.} Investment costs for existing technologies are assumed to be zero. For new technologies the investment costs are reported in Table \ref{table:CostAssumpstions_Technologies}. Fixed OPEX is calculated as a fraction of annualized CAPEX. Variable costs are given in terms of output of the respective technology. Offshore technologies are assumed to be 1.2x as expensive as their onshore counterpart.

\input{Tables/02_Technology_Cost_Assumptions.tex}

\textbf{Network Costs} Table \ref{table:CostAssumptions_Networks} shows cost assumptions for networks. DC links can only be installed at a rated power of 2 GW, in line with international standardization. The data has been recalculated for this study to fit the following formula, where $S$ is the capacity in MW and $d$ the distance in km.

\begin{equation}
    C_{c,l} = \gamma_1 + \gamma_2 S_\mathrm{l} + \gamma_3 \mathrm{d_l} + \gamma_4 \mathrm{d_l} S_\mathrm{l}
\end{equation}

\input{Tables/02_Network_Cost_Assumptions.tex}

\section{Energy System Model}
\subsection{Technology Performance}
The following technology models are used in this work:
\begin{itemize}
    \item \textit{Renewable, non-dispatchable technologies}: Electricity output that can be curtailed (Wind (onshore and offshore), Solar PV, Run of River, Biofuels)
    \item \textit{Conversion Technologies Type 1}: No input with electricity output. Fuel costs are allocated to the variable O\&M costs (Nuclear Power Plant, Coal \& Lignite Power Plant, Oil Power Plant)
    \item \textit{Conversion Technologies Type 2}: Input and output, fuel costs are separated from variable O\&M costs (Electrolyser, Fuel Cell, Gas Power Plant)
    \item \textit{Storage Technologies Type 1}: Most simple storage technology with input and output and an energy loss over time (Battery Storage, Hydro Storage - Closed Loop).
    \item \textit{Storage Technologies Type 2.1}: Same as Storage Technologies Type 1 but with an additional exogenous energy inflow over time, as for open-loop hydro technologies (Hydro Storage - Open Loop, Hydro Storage - Reservoir).
    \item \textit{Storage Technologies Type 2.2}: Same as Storage Technologies Type 1 but with an energy requirement when charging (e.g. for compression, Storage - Hydrogen).
\end{itemize}

Subsequently, we do not denote indices for nodes and technologies for readability.

\textbf{Renewable, non-dispatchable technologies}. The maximum possible output $\overline{\mathrm{X}}_\mathrm{out,el,t}$ of this technology can be curtailed. The maximum possible output is generated from the reneable generation profiles described in the last section As such the output in time-step $t$ is:
\begin{equation}
    \label{eq:technologies_RE}
    X_\mathrm{out,el,t} \leq \overline{\mathrm{X}}_\mathrm{out,el,t}
\end{equation}

\textbf{Conversion Technologies Type 1}.
The output of this technology type is only limited by its size. As such:
\begin{equation}
    \label{eq:technologies_type1}
    X_\mathrm{out,el,t} \leq S
\end{equation}

\textbf{Conversion Technologies Type 2}.
The output of this technology type is limited by its size and relates to the sum of its input. 
\begin{equation}
    \label{eq:technologies_type2_S}
    \sum X_\mathrm{out, r, t} \leq S
\end{equation}
\begin{equation}
    \label{eq:technologies_type2_InOut}
    X_\mathrm{out, r, t} = \alpha \sum X_\mathrm{in, r, t}
\end{equation}
Note that a natural gas power plant can have two inputs (hydrogen and natural gas), whereas the input of hydrogen is limited to $\kappa_\mathrm{r} = 5\%$. As such:
\begin{equation}
    \label{eq:technologies_type2_limIn}
    X_\mathrm{in, r, t} \leq \kappa_\mathrm{r} \sum X_\mathrm{in, r, t}
\end{equation}

\textbf{Storage Technologies Type 1}
All storage technologies have a maximum charging $\underline{\mathrm{x}}_\mathrm{in,r}$ and discharging rate $\overline{\mathrm{x}}_\mathrm{out,r}$ as well as a maximum charging level $S$:
\begin{equation}
    \label{eq:storage1_maxcharge}
    X_\mathrm{in, r, t} \leq \underline{\mathrm{x}}_\mathrm{in,r}S
\end{equation}
\begin{equation}
    \label{eq:storage1_maxdischarge}
    X_\mathrm{out, r, t} \leq \overline{\mathrm{x}}_\mathrm{out,r}S
\end{equation}
\begin{equation}
    \label{eq:storage1_size}
    S_\mathrm{r, t} \leq S
\end{equation}
 The state of charge $S_\mathrm{r, t}$ is connected to the charging and discharging of the storage as follows. Additionally a fraction of stored energy is lost over time at the share $\lambda$:
\begin{equation}
    \label{eq:storage1_stateofcharge}
       S_\mathrm{r, t} =  (1 - \lambda)S_\mathrm{r, t-1} + \eta_\mathrm{in} X_\mathrm{in, r, t} - 1 / \eta_\mathrm{out} X_\mathrm{out, r, t}
\end{equation}

To ensure that only either charging or discharging happens, we formulated an additional cut. While charging and discharging at the same time cannot be not completely avoided this way, it decreases the size of feasible solutions significantly and removes a corner point in the resulting MILP. For performance reasons we refrain from formulating a formulation involving binaries in this work. 
\begin{equation}
    \label{eq:storage1_cut}
    X_\mathrm{in, r, t} / \underline{\mathrm{x}}_\mathrm{in,r} + X_\mathrm{out, r, t} / \overline{\mathrm{x}}_\mathrm{out,r} \leq S_\mathrm{r, t}
\end{equation}

 \textbf{Storage Technologies Type 2.1}
 This technology model is used for hydro storage technologies that have a natural, exogenous inflow $\mathrm{X}_\mathrm{inflow, r, t}$ into the upper reservoir. As such, equation \ref{eq:storage1_stateofcharge} is adapted as follows:
 \begin{equation}
    \label{eq:storage1.1_stateofcharge}
       S_\mathrm{r, t} =  (1 - \lambda)S_\mathrm{r, t-1} + \eta_\mathrm{in} X_\mathrm{in, r, t} - 1 / \eta_\mathrm{out} X_\mathrm{out, r, t} + \mathrm{X}_\mathrm{inflow, r, t}
\end{equation}

 \textbf{Storage Technologies Type 2.2}
 This technology model is used for hydrogen storage that requires extra electricity for compression. As such, an additional constraint is added calculating the respective energy requirements:
\begin{equation}
    \label{eq:storage1.2_energyinput}
    X_\mathrm{in, el, t} = \gamma X_\mathrm{in, h2, t}
\end{equation}

A summary of all technology performance parameters is shown in the table below.

\input{Tables/02_Technology_Performance_Assumptions}

\subsection{Network Performance}
The constraints for the network performance are formulated for each branch respectively. We refrain from denoting indices on the branch the network for readability respectively.
The outflow from a node (i.e. the inflow to the branch) through a branch cannot exceed its capacity $S$:
\begin{equation}
    \label{eq:network_flow}
    F_\mathrm{out, n_1, n_2 , r, t} \leq S_\mathrm{n_1, n_2}
\end{equation}
At the other end of the branch we account for the losses through the branch. They are calculated based on a loss factor $\mu$ and the length of the branch $\mathrm{d}$.
\begin{equation}
    \label{eq:network_loss}
    F_\mathrm{in, n_1, n_2, r, t} = \mu \mathrm{d} F_\mathrm{out, n_1, n_2, r, t}
\end{equation}

For \textbf{all electricity networks}, bi-directional flows are allowed and sizes in both directions need to be equal.
\begin{equation}
    \label{eq:network_size_bidirectional}
    S_\mathrm{n_1, n_2} = S_\mathrm{n_2, n_1}
\end{equation}

To ensure only a one-directional flow per time-slice we follow a similar approach as for the storage technologies (eq. \ref{eq:storage1_cut}):
\begin{equation}
    \label{eq:network_cut_bidirectional}
    F_\mathrm{in, n_1, n_2, r, t} + F_\mathrm{in, n_2, n_1, r, t} \leq S_\mathrm{n_1, n_2}
\end{equation}

For \textbf{all hydrogen networks}, two branches need to be build for each direction respectively. Additionally, there is an electricity consumption for compression:
\begin{equation}
    \label{eq:network_consumption_unidirectional}
   F_\mathrm{cons, n_1, el, t} = \mathrm{k} F_\mathrm{out, n_1, n_2 , r, t}
\end{equation}

$\mathrm{k}$ is determined with:

\begin{equation}
    \label{eq:network_consumption_parameter}
   \mathrm{k} = \frac{\mathrm{c} \mathrm{T}}{\eta \mathrm{LHV}} \left(\frac{p}{30 \mathrm{bar}}\right)^{\frac{(\gamma-1)}{\gamma}-1}
\end{equation}

The table below shows the network performance parameters.

\input{Tables/02_Network_Performance}

\subsection{Managing computational complexities}
All optimizations have been run on full and hourly temporal resolution. For the 2030 cases, DC grids can only be build in 1GW capacities. To reduce the computational burden, this assumption has been relaxed to a continuous size for the 2040 cases. 

In the 2040 cases, we additionally used warm-starts for the \textit{Synergies} and the \textit{H-2 (only offshore)} scenario. We therefore used the solution of the best scenario (in this case \textit{T-All}). To manage numerical issues, we followed a three-step approach:

\begin{enumerate}
    \item Fix technology and network sizes to the solution from the \textit{T-All scenario}, fix all new technology and network sizes to zero and solve the problem. This converges to the same solution as the \textit{T-All} scenario, but provides a feasible solution for the next stage.
    \item Unfix variables of new technology and network sizes, but keep the sizes of technologies and networks in \textit{T-All} fixed. Use the solution and information from the last solution as a warm start.
    \item Unfix all variables and use the solution and information from the last solution as a warm start. This provides the optimal solution for the \textit{Synergies} scenario.
\end{enumerate}

\begin{landscape}
\section{Installed Capacities for the 2030 scenarios}
\input{Tables/02_ResultsInstalledCapacities2030}
\pagebreak
\section{Installed Capacities for the 2040 scenarios}
\input{Tables/02_ResultsInstalledCapacities2040}
\end{landscape}


\end{document}

%% file: Tables/00_InstalledCapacities_Countries.tex
\begin{table}
\caption{Considered electricity generation and storage technologies and networks in the starting system and their and capacities per country in GW.}
\label{tab:00_InstalledCapacities}
\resizebox{\textwidth}{!}{%
\begin{tabular}{lrrrrrrr}
\toprule
& BE & DE & DK & NL & NO & UK & total \\
\midrule
Dispatchable Technologies &  &  &  &  &  &  &  \\
\SmallIndent Coal \& Lignite & - & 17.0 & - & - & - & - & 17.0 \\
\SmallIndent Gas & 10.1 & 35.8 & 1.7 & 14.4 & 0.3 & 27.2 & 89.5 \\
\SmallIndent Nuclear & 2.1 & - & - & 0.5 & - & 9.3 & 11.9 \\
\SmallIndent Oil & - & 0.9 & - & - & - & - & 0.9 \\
\SmallIndent Hydro Storage - Closed Loop (Energy Cap.)* & 5.8 & 629.9 & - & - & - & 26.4 & 662.0 \\
\SmallIndent Hydro Storage - Closed Loop (Pump Cap.) & 1.2 & 7.4 & - & - & - & 2.7 & 11.3 \\
\SmallIndent Hydro Storage - Closed Loop (Turbine Cap.) & 1.3 & 7.4 & - & - & - & 2.7 & 11.4 \\
\SmallIndent Hydro Storage - Open Loop (Energy Cap.)* & - & 416.8 & - & - & 89577.4 & - & 89994.2 \\
\SmallIndent Hydro Storage - Open Loop (Pump Cap.) & - & 1.4 & - & - & 1.1 & - & 2.5 \\
\SmallIndent Hydro Storage - Open Loop (Turbine Cap.) & - & 1.6 & - & - & 37.8 & - & 39.5 \\
\SmallIndent Hydro Storage - Reservoir (Energy Cap.)* & - & 258.6 & - & - & - & - & 258.6 \\
\SmallIndent Hydro Storage - Reservoir (Turbine Cap.) & - & 1.3 & - & - & - & - & 1.3 \\
\midrule
Non-dispatchable Technologies &  &  &  &  &  &  &  \\
\SmallIndent Wind Offshore & 5.9 & 30.5 & 6.9 & 15.5 & - & 48.9 & 107.7 \\
\SmallIndent Wind Onshore & 4.7 & 75.4 & 6.2 & 8.0 & 6.1 & 26.6 & 126.8 \\
\SmallIndent Solar & 10.4 & 96.1 & 6.5 & 27.3 & 0.6 & 23.4 & 164.3 \\
\SmallIndent Run of River (Turbine Cap.) & 0.1 & 4.4 & - & - & - & 2.1 & 6.6 \\
\SmallIndent Biofuels & 0.9 & 7.6 & 2.2 & 1.2 & - & 6.4 & 18.4 \\
\midrule
Networks &  &  &  &  &  &  &  \\
\SmallIndent Electricity (AC) & \multicolumn{7}{c}{see Figure \ref{fig:SystemTopology}} \\
\SmallIndent Electricity (DC) &  &  &  &  &  &  &  \\
\bottomrule
\multicolumn{8}{l}{* Energy Capacity in GWh}
\end{tabular}
}
\end{table}

%% file: Tables/00_DataSources.tex
\renewcommand{\arraystretch}{1.4}
\begin{table}[]
\caption{Summary of data collected for this work with respective sources}
\label{tab:DataSources}
\begin{tabularx}{\textwidth}[t]{p{3.6cm} X p{1.5cm}}
\toprule
 \textbf{Data Set} & \textbf{Description} & \textbf{References}\\
\midrule
\multicolumn{3}{l}{\textbf{Demand}}\\
Electricity & Based on hourly, national profiles from \acs{TYNDP} National Trend scenario for the climate year of 2008. Allocation to nodes based on spatial distribution of total annual demand and industrialization. & \cite{ENTSO-E2022a, Neumann2023, Eurostat2023a,Eurostat2023, Eurostat2023b}\\
Hydrogen & Based on annual demand on NUTS2 level and distributed equally as a flat profile over the year. & \cite{Gross2022}\\
\midrule
\multicolumn{3}{l}{\textbf{Installed capacities - Conversion and storage technologies}} \\
\parbox[t]{\textwidth}{Coal \& lignite,\\ Gas,\\ Nuclear,\\ Oil,\\ Hydro storage - closed loop,\\ Hydro storage - open loop,\\ Hydro storage - reservoir} & National capacities are based on projected capacities in 2030 as reported in \acs{TYNDP} 2022 for the scenario National Trend. Allocation to nodes using node keys derived from a spatially resolved power plant database. &  \cite{ENTSO-E2022a, Gotzens2019} \\
Wind offshore & National data on wind park tenders, expected capacities and date of operation was used to form wind farm clusters. & National sources, \cite{4COffshore2023,NorthSeaEnergy2023,EuropeanCommission2023}\\
\parbox[t]{\textwidth}{Wind onshore, \\ Solar PV} & Installed capacities for 2023 and potentials on a NUTS2 level are used to project installed capacities in 2030 per node using projected capacities in the \acs{TYNDP} National Trend scenario. & National sources, \cite{ENTSO-E2022a, JointResearchCentre2021} \\
Run of river & National capacities are based on projected capacities in 2030 as reported in \acs{TYNDP} 2022 for the scenario National Trend. Allocation to nodes using node keys derived from a spatially resolved power plant database. &  \cite{ENTSO-E2022a, Gotzens2019} \\
Biofuels & Installed capacities for 2023 and potentials on a NUTS2 level are used to project installed capacities in 2030 per node using projected capacities in the \acs{TYNDP} National Trend scenario. & National sources, \cite{ENTSO-E2022a, JointResearchCentre2021} \\
\midrule
\multicolumn{3}{l}{\textbf{Installed capacities - Networks}}\\
\parbox[t]{\textwidth}{Electricity (AC), \\ Electricity (DC)} & Based on  grid data used in Neumann et al. (2023) and corrected using data from ENTSO-E and Hoogspanningsnet. Offshore transmission capacities collected from national sources. & National sources, \cite{Neumann2023, Hoogspanningsnet2024, ENTSO-E2022}\\
\midrule
\multicolumn{3}{l}{\textbf{Renewable generation profiles}}\\
Wind offshore & Centroids of each wind farm, wind speed data from ERA5 and turbine types for each wind farm was used. & \cite{Munoz-Sabater2021}\\
Wind onshore & Installed capacities for 2030 on NUTS2 level and climate data from ERA5 for 2008 was used to calculate hourly generation profiles. Power curve for an average 1.5MW turbine with a hub height of 100m was used. & \cite{Munoz-Sabater2021, ENTSO-E2022a, ENTSO-E2022}\\
Solar PV & Installed capacities for 2030 on NUTS2 level and climate data from ERA5 for 2008 was used to calculate hourly generation profiles & \cite{Munoz-Sabater2021}\\
Run of river & Based on daily generation from  ERAA and capacity per node. The daily generation was equally distributed to the hours of each day. & \cite{ENTSO-E2022}\\
Biofuels & A capacity factor at each node of 0.53 was assumed. & \cite{Bolson2022}\\
\midrule
\multicolumn{3}{l}{\textbf{Natural inflows to hydro storage technologies}}\\
\parbox[t]{\textwidth}{Hydro storage - open loop, \\Hydro storage - reservoir}  & Weekly inflow values from ERAA were equally distributed over the hours of a week. Allocation of national inflows to each node is based on the capacity installed per node & \cite{ENTSO-E2022}\\
\midrule
\multicolumn{3}{l}{\textbf{Network expansion limits}}\\
Electricity & Based on  grid data used in Neumann et al. (2023) and corrected using data from ENTSO-E and Hoogspanningsnet. Offshore transmission capacities collected from national sources. & National sources, \cite{Neumann2023, Hoogspanningsnet2024, ENTSO-E2022}\\
Hydrogen & Based on data from publications of the European Hydrogen Backbone and national sources. Existing onshore and offshore pipelines are taken into account. & National Sources, \cite{EHB2022}\\
\midrule
\multicolumn{3}{l}{\textbf{Cost assumptions}}\\
Technologies & Based on data from the Danish energy agency, and Neumann et al (2022) & \cite{DanishEnergyAgency2022, Neumann2023}\\
Networks & Based on data from Zappa et al (2019) and EHB (2022) & \cite{Zappa2019, EHB2022}\\
\bottomrule
\end{tabularx}
\end{table}

%% file: Tables/00_ScenarioDefinitions.tex
\renewcommand{\arraystretch}{1}


\begin{table}[]
\caption{Scenario Definitions}
\label{tab:ScenarioDefinition}
\resizebox{\textwidth}{!}{%
\begin{tabular}{ll}
  \toprule
   & 
   \begin{tabularx}{\textwidth}{L{1.5cm}X*{3}{p{0.3cm}} p{0.15cm} *{2}{p{0.3cm}} p{0.15cm} *{5}{p{0.3cm}}}
    & & \multicolumn{3}{c}{\makecell{Grid\\Expansion}} & & \multicolumn{2}{c}{\makecell{Storage}} & & \multicolumn{5}{c}{\makecell{Hydrogen\\Technologies}} \\ 
    \cmidrule(l){3-5} \cmidrule(l){7-8} \cmidrule(l){10-14}
    Scenario & Description & \rotatebox{90}{Onshore} & \rotatebox{90}{Offshore} & \rotatebox{90}{Over borders} & & \rotatebox{90}{Onshore} & \rotatebox{90}{Offshore} & & \rotatebox{90}{Electrolysis (onshore)} & \rotatebox{90}{Electrolysis (offshore)} & \rotatebox{90}{H2 Storage} & \rotatebox{90}{Fuel Cell} & \rotatebox{90}{Pipelines} \\
   \end{tabularx}\\
   \midrule
   \rotatebox[origin=c]{90}{} & 
   \begin{tabularx}{\textwidth}{L{1.5cm}X*{3}{p{0.3cm}} p{0.15cm} *{2}{p{0.3cm}} p{0.15cm} *{5}{p{0.3cm}}}
       Reference & 2030: No additional technologies or networks possible. Resembles roughly the energy system in 2030 as projected by TYNDP in the National Trends scenario. & & & & & & & & & & & & \\
       & 2040: Same as for 2030 case with increased electricity and hydrogen demand. Additionally, onshore wind, offshore wind and onshore PV capacities can be expanded & & & & & & & & & & & & \\

   \end{tabularx}\\
   %
   \midrule
   \rotatebox[origin=c]{90}{Role of transmission} & 
   \begin{tabularx}{\textwidth}{L{1.5cm}X*{3}{p{0.3cm}} p{0.15cm} *{2}{p{0.3cm}} p{0.15cm} *{5}{p{0.3cm}}}
      T-All & All transmission corridors can be used and expanded up to a given limit. New corridors are also available. & x & x & x & & & & & & & & & \\
      \midrule
      T1 & Only onshore corridors can be used and expanded. & x & & x & & & & & & & & & \\
      \midrule
      T2 & Only offshore corridors can be used and expanded. & & x & x & & & & & & & & & \\
      \midrule
      T3 & All transmission corridors can be used and expanded up to a given limit, except for corridors crossing a  national border. & x & x & & & & & & & & & & \\
   \end{tabularx}\\
  %
   \midrule
   \rotatebox[origin=c]{90}{Role of electricity storage} & 
   \begin{tabularx}{\textwidth}{L{1.5cm}X*{3}{p{0.3cm}} p{0.15cm} *{2}{p{0.3cm}} p{0.15cm} *{5}{p{0.3cm}}}
      S-All & Lithium-ion batteries can be installed at all nodes (onshore and offshore). At offshore nodes, there is a  limit of two large offshore platforms with a respective storage limit. & & & & & x & x & & & & & \\
      \midrule
      S1 & Only onshore nodes are available for energy storage & & & & & x & & & & & & & \\
      \midrule
      S2 & Only offshore nodes are available for energy storage & & & & & & x & & & & & & \\
      \midrule
      S-All-HPE & Same as S-All, but with a power-to-energy ratio of 1, i.e. all energy stored can be discharged within one hour. & & & & & x & x & & & & & & \\
   \end{tabularx}\\
   %
   \midrule
   \rotatebox[origin=c]{90}{Role of hydrogen} & 
   \begin{tabularx}{\textwidth}{L{1.5cm}X*{3}{p{0.3cm}} p{0.15cm} *{2}{p{0.3cm}} p{0.15cm} *{5}{p{0.3cm}}}
      H-All & All hydrogen technologies can be installed. This includes production, storage, transport and reconversion into electricity. & & & & & & & & x & x & x & x & x \\
      \midrule
      H2 & Hydrogen production and transport only onshore & & & & & & & & x & & x & x & x \\
      \midrule
      H3 & Hydrogen production only offshore & & & & & & & & & x & x & x & x \\
      \midrule
      H3 & Excludes hydrogen storage & & & & & & & & x & x & & x & x \\
      \midrule
      H4 & Excludes hydrogen transport: hydrogen can only be used at the node where it is produced. As such, also offshore hydrogen production is excluded. & & & & & & & & x & & x & x & \\
   \end{tabularx}\\
   \midrule
   \rotatebox[origin=c]{90}{} & 
   \begin{tabularx}{\textwidth}{L{1.5cm}X*{3}{p{0.3cm}} p{0.15cm} *{2}{p{0.3cm}} p{0.15cm} *{5}{p{0.3cm}}}
      Synergies & All technologies and networks of previous scenarios can be expanded/newly build & x & x & x & & x & x & & x & x & x & x & x \\
   \end{tabularx}\\
   \bottomrule
\end{tabular}
}
\end{table}

%% file: Tables/01_ResultsReference.tex
\begin{table}[h]
\caption{Results of the \textit{Reference} scenario}
\label{tab:Results_Reference}
\centering
\resizebox{0.8\textwidth}{!}{%
\begin{tabular}{lrrrrrc}
\toprule
Country & \multicolumn{3}{c}{Carbon Emissions (Mt)} & \multicolumn{2}{c}{Renewable Share} & \\
  & From H2 Production & From Electricity Generation & Total & Results & National Goal & Source \\
 \midrule
BE & 2.50 & 6.26 & 8.76 & 59.1\% & 37.4\% & \cite{IEAInternationalEnergyAgency2022} \\
DE & 13.44 & 36.02 & 49.46 & 79.6\% & 80\% & \cite{Bundesregierung2023} \\
DK & 0.58 & 0.55 & 1.13 & 97.3\% & 117\%$^\mathrm{a}$ & \cite{Energistyrelsen2023} \\
NL & 4.63 & 8.05 & 12.68 & 79.1\% & 70\% & \cite{MinistryofEconomicAffairsandClimatePolicy2019} \\
NO & 0.46 & 0.00 & 0.46 & 100.0\% & 100\% & \cite{IEAInternationalEnergyAgency2022a} \\
UK & 7.86 & 6.35 & 14.2 & 84.2\% & 95\%$^\mathrm{b}$ & \cite{IEAInternationalEnergyAgency2022a} \\
\bottomrule
\multicolumn{3}{c}{\makecell[l]{$^\mathrm{a}$ includes cross border flows\\ $^\mathrm{b}$ includes power generation from nuclear power plants}}
\end{tabular}
}
\end{table}

%% file: Tables/02_InstalledCapacities_DataSources.tex
\begin{table}[h]
\centering
\caption{Sources for installed capacities in 2030 and their allocation}
\label{tab:02_InstalledCapacities_DataSources}
\begin{tabular}{lP{2.5cm}P{2.5cm}P{4.5cm}}
\toprule
Technology & Source national capacity & Allocation Keys per Node & Note \\
\midrule
Oil & \cite{ENTSO-E2022a} & \cite{Gotzens2019}\\
Gas & \cite{ENTSO-E2022a} & \cite{Gotzens2019}\\
Nuclear & \cite{ENTSO-E2022a} & \cite{Gotzens2019}\\
Coal \& Lignite & \cite{ENTSO-E2022a} & \cite{Gotzens2019}\\
Hydro (run of river) & \cite{ENTSO-E2022a} & \cite{Gotzens2019} & Treated as non dispatchable\\
Hydro (closed loop) &  \cite{ENTSO-E2022} & \cite{Gotzens2019} & Allocation key for size, turbine and pump power the same\\
Hydro (open loop) &  \cite{ENTSO-E2022} & \cite{Gotzens2019} & Allocation key for size, turbine and pump power the same\\
Hydro (reservoir) &  \cite{ENTSO-E2022} & \cite{Gotzens2019} & Allocation key for size, turbine and pump power the same\\
PV & \cite{ENTSO-E2022a} & National Sources & Treated as non dispatchable\\
Wind onshore & \cite{ENTSO-E2022a} & National Sources & Treated as non dispatchable\\
Wind offshore & Treated as non dispatchable\\
Other non-dispatchable RE & \cite{ENTSO-E2022a} & National Sources & Treated as non dispatchable\\
\bottomrule
\end{tabular}
\end{table}

%% file: Tables/02_InstalledCapacities_NonOffshore_Country.tex
\begin{table}
\centering
\caption{Installed Capacities in GW (aggregated per Country and per source)}
\begin{tabular}{llrrr}
\toprule
   &              &  Capacity our work &  Capacity ENTSO-E &  Capacity PyPsa \\
Country & Technology &                    &                   &                 \\
\midrule
BE & Biofuels &               0.90 &              0.90 &            0.00 \\
   & Gas &              10.11 &             10.11 &            5.34 \\
   & Hydro - Pump Storage Closed Loop (Energy) &               5.75 &              5.75 &            0.00 \\
   & Hydro - Pump Storage Closed Loop (Turbine) &               1.30 &              1.30 &            1.31 \\
   & Hydro - Run of River (Turbine) &               0.15 &              0.15 &            0.05 \\
   & Nuclear &               2.08 &              2.08 &            2.08 \\
   & PV &              10.40 &             10.40 &            0.00 \\
   & Wind Onshore &               4.67 &              4.67 &            0.00 \\
DE & Biofuels &               7.57 &              7.57 &            0.86 \\
   & Coal \& Lignite &              17.04 &             17.04 &           26.20 \\
   & Gas &              35.77 &             35.77 &           25.29 \\
   & Hydro - Pump Storage Closed Loop (Energy) &             629.88 &            629.88 &            0.00 \\
   & Hydro - Pump Storage Closed Loop (Turbine) &               7.38 &              7.38 &            7.59 \\
   & Hydro - Pump Storage Open Loop (Energy) &             416.76 &            416.76 &            0.00 \\
   & Hydro - Pump Storage Open Loop (Turbine) &               1.64 &              1.64 &            0.00 \\
   & Hydro - Reservoir (Energy) &             258.58 &            258.58 &            0.00 \\
   & Hydro - Reservoir (Turbine) &               1.30 &              1.30 &            0.17 \\
   & Hydro - Run of River (Turbine) &               4.37 &              4.37 &            2.05 \\
   & Oil &               0.85 &              0.85 &            0.84 \\
   & PV &              96.14 &             96.14 &            0.00 \\
   & Wind Onshore &              75.37 &             75.37 &            0.00 \\
DK & Biofuels &               2.24 &              2.24 &            0.00 \\
   & Gas &               1.66 &              1.66 &            1.70 \\
   & PV &               6.47 &              6.47 &            0.00 \\
   & Wind Onshore &               6.16 &              6.16 &            0.00 \\
NL & Biofuels &               1.25 &              1.25 &            0.00 \\
   & Gas &              14.43 &             14.43 &            9.10 \\
   & Nuclear &               0.49 &              0.49 &            0.48 \\
   & PV &              27.26 &             27.26 &            0.00 \\
   & Wind Onshore &               7.98 &              7.98 &            0.00 \\
NO & Gas &               0.27 &              0.27 &            0.64 \\
   & Hydro - Pump Storage Open Loop (Energy) &           89577.42 &          89577.42 &            0.00 \\
   & Hydro - Pump Storage Open Loop (Turbine) &              37.83 &             37.83 &            0.00 \\
   & PV &               0.60 &              0.60 &            0.00 \\
   & Wind Onshore &               6.06 &              6.06 &            0.00 \\
UK & Biofuels &               6.44 &              6.44 &            0.39 \\
   & Gas &              27.24 &             27.24 &           32.64 \\
   & Hydro - Pump Storage Closed Loop (Energy) &              26.38 &             26.38 &            0.00 \\
   & Hydro - Pump Storage Closed Loop (Turbine) &               2.74 &              2.74 &            0.44 \\
   & Hydro - Reservoir (Turbine) &               0.01 &              0.01 &            0.13 \\
   & Hydro - Run of River (Turbine) &               2.10 &              2.10 &            8.87 \\
   & Nuclear &               9.33 &              9.33 &            1.32 \\
   & PV &              23.41 &             23.41 &            0.00 \\
   & Wind Onshore &              26.59 &             26.59 &            0.00 \\
\bottomrule
\end{tabular}
\end{table}

%% file: Tables/02_InstalledCapacities_BiomassAllocation.tex
\begin{table}
    \centering
    \caption{Allocation of biomass capacities}
    \label{tab:biomass_allocation}
    \begin{tabular}{lll}
        \toprule
         Country & Note & Sources\\ 
         \midrule
         United Kingdom & Capacity allocation according to source & \cite{UKGovernment2023} \\
         Belgium & All capacity allocated to BE2 & \\
         The Netherlands & Capacity of Power Plant Eemshaven allocated to NL3. Rest equally distributed & \\
         Germany & Capacity allocation according to source & \cite{Bundesnetzagentur2023} \\
         Denmark & No allocation required (one node only) & \\
         Norway & No allocation required (one node only) & \\
      \bottomrule
    \end{tabular}
\end{table}

%% file: Tables/02_InstalledCapacities_OffshoreWindSources.tex
\begin{table}
    \centering
    \caption{Sources for offshore wind farms}
    \label{tab:offshore_wind_sources}
    \begin{tabular}{cc}
        \toprule
         Country & Sources\\ 
         \midrule
         United Kingdom & \cite{TheCrownEstate,TheCrownEstate2023,TheCrownEstateScotland2023} \\ 
         Belgium & \cite{Economie.fgov2023}\\
         Netherlands & \cite{MinisterievanEconomischeZakenenKlimaat2022,GovernmentoftheNetherlands2020,RVO2023} \\
         Germany & \cite{BundesamtfurSeeschifffahrtundHydrographie2023,NetzentwicklungsplanStrom2023} \\
         Denmark & \cite{Dataforsyningen,ENS2023} \\
         Norway & \cite{NVE2023a,NVE2023} \\
         \bottomrule
    \end{tabular}
\end{table}

%% file: Tables/02_Demand_Country.tex
\begin{table}[h]
\centering
\caption{Annual demand in TWh (aggregated per Country for the year 2030 and climate year 2009)}
\begin{tabular}{lrrr}
\toprule
{} &  Our Work/TYNDP &  Neumann et al. (2023) \cite{Neumann2023} &  Eurostat (2023), demand in 2019 \cite{Eurostat2023c} \\
\midrule
BE &    95 &                 131 &             82 \\
DE &   585 &                 730 &            497 \\
DK &    53 &                  50 &             31 \\
NL &   140 &                 186 &            108 \\
NO &   165 &                 113 &            116 \\
\bottomrule
\end{tabular}
\end{table}

%% file: Tables/02_DemandSupply_country.tex
\begin{table}
\centering
\caption{Capacity factors, total annual generation from non-dispatchable sources without curtailment and total demand per country in TWh}
\label{tab:02_DemandSupplyCountry}
\begin{tabular}{llrrrrrr}
\toprule
      & Country &    BE &     DE &    DK &     NL &     NO &     UK \\
Technology & Type &       &        &       &        &        &        \\
\midrule
Biomass & CF &  0.53 &   0.53 &  0.53 &   0.53 &   0.53 &   0.53 \\
      & Generation &  4.20 &  35.15 & 10.40 &   5.80 &   0.00 &  29.90 \\
PV & CF &  0.10 &   0.10 &  0.10 &   0.09 &   0.11 &   0.09 \\
      & Generation &  8.70 &  82.36 &  5.46 &  22.46 &   0.59 &  18.27 \\
Run of River & CF &  0.32 &   0.53 &     0 &      0 &      0 &   0.24 \\
      & Generation &  0.42 &  20.17 &  0.00 &   0.00 &   0.00 &   4.50 \\
Wind offshore & CF &  0.50 &   0.58 &  0.53 &   0.53 &      0 &   0.55 \\
      & Generation & 25.96 & 153.92 & 31.65 &  71.91 &   0.00 & 235.78 \\
Wind onshore & CF &  0.21 &   0.16 &  0.30 &   0.29 &   0.38 &   0.23 \\
      & Generation &  8.75 & 106.16 & 16.01 &  20.25 &  19.96 &  52.89 \\
\midrule
Total Generation &  & 48.02 & 397.77 & 63.53 & 120.41 &  20.55 & 341.34 \\
Demand &  & 95.45 & 585.09 & 53.25 & 139.65 & 165.41 & 306.70 \\
\bottomrule
\end{tabular}
\end{table}

%% file: Tables/02_Technology_Cost_Assumptions.tex
\begin{table}[h]
\centering
\centering
\caption{Technology Cost assumptions.}
\label{table:CostAssumpstions_Technologies}
\resizebox{\textwidth}{!}{%
\begin{tabular}{lrrrrrl}
	\toprule
	\multicolumn{1}{l}{\makecell{Technology}} & \multicolumn{1}{c}{\makecell{CAPEX}} & \multicolumn{1}{c}{\makecell{Lifetime}} & \multicolumn{1}{c}{\makecell{OPEX Fixed}} & \multicolumn{1}{c}{\makecell{OPEX Variable}} & \multicolumn{1}{c}{\makecell{Emission Factor}} & \multicolumn{1}{c}{\makecell{Source}}\\
 	\multicolumn{1}{c}{\makecell{}} & \multicolumn{1}{c}{\makecell{[kEUR]}} & \multicolumn{1}{c}{\makecell{[a]}} & \multicolumn{1}{c}{\makecell{[\% annual CAPEX]}} & \multicolumn{1}{c}{\makecell{[EUR]}} & \multicolumn{1}{c}{\makecell{[kg/MWh$_\mathrm{el}$]}} & \multicolumn{1}{c}{\makecell{}}\\
	\midrule
 Existing Technologies & & & & & &\\
\SmallIndent Power Plant (Coal) & & & & 62.84$^4$ & 763 & \cite{DanishEnergyAgency2022, Neumann2023}\\
\SmallIndent Power Plant (Gas) & & & & 4.20 & 302 & \cite{DanishEnergyAgency2022, Neumann2023}\\
\SmallIndent Power Plant (Oil) & & & & 130.23$^4$ & 631 & \cite{Neumann2023}\\
\SmallIndent Power Plant (Nuclear) & & & & 16.90$^4$ & 0 & \cite{Neumann2023} \\
\SmallIndent Pumped Hydro (Closed) & & & & 0.00 & 0 & \cite{DanishEnergyAgency2022} \\
\SmallIndent Pumped Hydro (Open) & & & & 0.00 & 0 & \cite{DanishEnergyAgency2022} \\
\SmallIndent Pumped Hydro (Reservoir) & & & & 0.00 & 0 & \cite{DanishEnergyAgency2022} \\
\midrule
 New Technologies & & & & & &\\
\SmallIndent Electrolyser (PEM) & 650$^1$ & 25 & 2.00\% & 0.00 & 0 & \cite{DanishEnergyAgency2022} \\
\SmallIndent Electrolyser (PEM, Offshore) & 780$^1$ & 25 & 2.00\% & 0.00 & 0 & \cite{DanishEnergyAgency2022} \\
\SmallIndent Electrolyser (Alkaline) & 450$^1$ & 30 & 4.00\% & 0.00 & 0 & \cite{DanishEnergyAgency2022} \\
\SmallIndent Battery System (Onshore) & 622$^2$ & 25 & 7.91\% & 1.80 & 0 & \cite{DanishEnergyAgency2022} \\
\SmallIndent Battery System (Offshore) & 746$^2$ & 25 & 7.91\% & 1.80 & 0 & \cite{DanishEnergyAgency2022} \\
\SmallIndent Fuel Cell & 1100$^3$ & 10 & 5.00\% & 0.00 & 0 & \cite{DanishEnergyAgency2022} \\
\SmallIndent Hydrogen Storage (Cavern) & 2$^2$ & 100 & 0.00\% & 0.00 & 0 & \cite{DanishEnergyAgency2022} \\
\bottomrule
\multicolumn{7}{l}{\parbox{\linewidth}{$^1$ Based on MW$_\mathrm{el}$ input}}\\
\multicolumn{7}{l}{\parbox{\linewidth}{$^2$ Based on MWh storage capacity}}\\
\multicolumn{7}{l}{\parbox{\linewidth}{$^3$ Based on MW$_\mathrm{el}$ output}}\\
\multicolumn{7}{l}{\parbox{\linewidth}{$^4$ Including fuel costs}}\\
\end{tabular}
}
\end{table}

%% file: Tables/02_Network_Cost_Assumptions.tex
\begin{table}[h]
\centering
\centering
\caption{Network Cost assumptions.}
\label{table:CostAssumptions_Networks}
\resizebox{\textwidth}{!}{%
\begin{tabular}{lSSSSSSSSS}
\toprule
     & \multicolumn{1}{c}{\makecell{$\gamma_1$}} & \multicolumn{1}{c}{\makecell{$\gamma_2$}} & \multicolumn{1}{c}{\makecell{$\gamma_3$}} & \multicolumn{1}{c}{\makecell{$\gamma_4$}} & \multicolumn{1}{c}{\makecell{Lifetime}} & \multicolumn{1}{c}{\makecell{OPEX Fix}} & \multicolumn{1}{c}{\makecell{OPEX Var}} & \multicolumn{1}{c}{\makecell{Source}}\\
     & \multicolumn{1}{c}{\makecell{[kEUR]}} & \multicolumn{1}{c}{\makecell{[kEUR/ \\ MW]}} & \multicolumn{1}{c}{\makecell{[kEUR/ \\ km]}} & \multicolumn{1}{c}{\makecell{[kEUR \\ /km/MW]}} & \multicolumn{1}{c}{\makecell{[a]}} & \multicolumn{1}{c}{\makecell{[\% annual \\ CAPEX]}} & \multicolumn{1}{c}{\makecell{[EUR \\ /MWh]}} & \multicolumn{1}{c}{\makecell{}}\\
\midrule
    Electricity Networks & & & & & & & &\\
    \SmallIndent AC & 0 & 43.7 & 0 & 0.4 & 40 & 0.04 & 0 & \cite{Zappa2019}\\
    \SmallIndent DC & 0 & 68.1 & 0 & 0.1 & 40 & 0.04 & 0 & \cite{Zappa2019}\\
\midrule
    Hydrogen Pipelines & & & & & & & &\\
    \SmallIndent Offshore & 337045.7 & -33.1 & 0 & 0.5 & 50 & 0.04 & 0 & \cite{EHB2022}\\
    \SmallIndent Onshore, new & 198262.2 & -19.5 & 0 & 0.3 & 50 & 0.04 & 0 & \cite{EHB2022}\\
    \SmallIndent Onshore, repurposed & 39567.9 & -3.9 & 0 & 0.1 & 50 & 0.04 & 0 & \cite{EHB2022}\\
\bottomrule
\end{tabular}
}
\end{table}

%% file: Tables/02_Technology_Performance_Assumptions.tex
\begin{table}[h]
\centering
\centering
\caption{Technology Performance parameters.}
\label{table:PerformanceAssumpstions_Technologies}
\begin{tabular}{lSSSSSSS}
\toprule
\multicolumn{1}{l}{\makecell{Technology}} & \multicolumn{1}{c}{\makecell{$\alpha$
}} & \multicolumn{1}{c}{\makecell{$\underline{\mathrm{x}}_\mathrm{in,r}$
}} & \multicolumn{1}{c}{\makecell{$\overline{\mathrm{x}}_\mathrm{out,r}$
}} & \multicolumn{1}{c}{\makecell{$\lambda$
}} & \multicolumn{1}{c}{\makecell{$\eta_\mathrm{in}$
}} & \multicolumn{1}{c}{\makecell{$\eta_\mathrm{out}$
}} & \multicolumn{1}{c}{\makecell{$\gamma$}}\\
	\midrule
Conversion Technology Type 2 &  &  &  &  &  &  & \\
\SmallIndent Electrolyzer & 0.655 &  &  &  &  &  & \\
\SmallIndent Fuel Cell & 0.500 &  &  &  &  &  & \\
\SmallIndent Power Plant Gas & 0.610 &  &  &  &  &  & \\
\midrule
Storage Technology Type 1 &  &  &  &  &  &  & \\
\SmallIndent Battery Storage &  & 0.333 & 0.333 & 4.168E-05 & 0.985 & 0.975 & \\
\midrule
Storage Technology Type 2.1 &  &  &  &  &  &  & \\
\SmallIndent Hydrogen Storage &  & 0.333 & 0.333 & 0.000 & 0.990 & 1.000 & 0.008\\
\midrule
Storage Technology Type 2.2 &  &  &  &  &  &  & \\
\SmallIndent Pumped Hydro Storage &  & * & * & 0.000 & 0.890 & 0.890 & \\
\bottomrule
\multicolumn{8}{l}{\parbox{\linewidth}{$^*$ Depend on the installed capacity at the node}}\\
\end{tabular}
\end{table}

%% file: Tables/02_Network_Performance.tex
\begin{table}[h]
\centering
\centering
\caption{Network Performance parameters.}
\label{table:PerformanceAssumpstions_Networks}
\begin{tabular}{lSSSSSSS}
\toprule
\multicolumn{1}{l}{\makecell{Network}} & \multicolumn{1}{c}{\makecell{$\mu$
}} & \multicolumn{1}{c}{\makecell{$\mathrm{p}$
}} & \multicolumn{1}{c}{\makecell{$\mathrm{c}$
}} & \multicolumn{1}{c}{\makecell{$\mathrm{T}$
}} & \multicolumn{1}{c}{\makecell{$\eta$
}} & \multicolumn{1}{c}{\makecell{$\gamma$
}} & \multicolumn{1}{c}{\makecell{LHV
}}\\
\midrule
Electricity Networks &  &  &  &  &  &  & \\
\SmallIndent Electricity AC & 7e-05 &  &  &  &  &  & \\
\SmallIndent Electricity DC & 4e-05 &  &  &  &  &  & \\
\midrule
Hydrogen Pipelines &  &  &  &  &  &  & \\
\SmallIndent All hydrogen pipelines & 4e-05 & 140 & 0.00398 & 300 & 0.65 & 1.405 & 33.32\\
\bottomrule
\end{tabular}
\end{table}

%% file: Tables/02_ResultsInstalledCapacities2030.tex
\begin{table}[h!]
\caption{Installed Capacities for all scenarios aggregated for the whole energy system design in 2030. Existing assets are not reported in this table.}
\label{tab:00_InstalledCapacities2030}
\begin{tabular}{l*{11}{S}}
\toprule
  & \multicolumn{3}{c}{\ce{H2} pipeline} & \multicolumn{2}{c}{Electricity grid (AC)} & \multicolumn{2}{c}{Battery} & \multicolumn{2}{c}{Electrolyser} & \multicolumn{1}{c}{Fuel Cell} & \multicolumn{1}{c}{\ce{H2} storage}\\
 & \multicolumn{1}{c}{offshore} & \multicolumn{1}{c}{onshore (new)} & \multicolumn{1}{c}{onshore (repurposed)} & \multicolumn{1}{c}{AC} & \multicolumn{1}{c}{DC} & \multicolumn{1}{c}{offshore} & \multicolumn{1}{c}{onshore} & \multicolumn{1}{c}{offshore} & \multicolumn{1}{c}{onshore} & \multicolumn{1}{c}{} & \multicolumn{1}{c}{}\\
 & \multicolumn{1}{c}{(GWkm)} & \multicolumn{1}{c}{(GWkm)} & \multicolumn{1}{c}{(GWkm)} & \multicolumn{1}{c}{(GWkm)} & \multicolumn{1}{c}{(GWkm)} & \multicolumn{1}{c}{(GW)} & \multicolumn{1}{c}{(GW)} & \multicolumn{1}{c}{(GWh)} & \multicolumn{1}{c}{(GW)} & \multicolumn{1}{c}{(GW)} & \multicolumn{1}{c}{(GWh)}\\
\midrule
T-1 (only onshore) &  &  &  & 997.37 & 0.00 &  &  &  &  &  & \\
T-2 (only offshore) &  &  &  & 0.00 & 24940.13 &  &  &  &  &  & \\
T-3 (no border cross) &  &  &  & 444.46 & 0.00 &  &  &  &  &  & \\
T-All &  &  &  & 571.87 & 22357.57 &  &  &  &  &  & \\
\midrule
S-1 (only onshore) &  &  &  &  &  &  & 0.00 &  &  &  & \\
S-2 (only offshore) &  &  &  &  &  & 0.00 &  &  &  &  & \\
S-All &  &  &  &  &  & 0.00 & 0.00 &  &  &  & \\
\midrule
H-1 (only onshore) & 0.00 & 0.00 & 1206.26 &  &  &  &  &  & 15.03 & 0.00 & 159.01\\
H-2 (only offshore) & 1676.64 & 0.00 & 2571.35 &  &  &  &  & 13.29 &  & 0.00 & 37.35\\
H-3 (no storage) & 0.00 & 0.00 & 1151.05 &  &  &  &  & 0.00 & 13.42 & 0.00 & \\
H-4 (only local use) &  &  &  &  &  &  &  &  & 14.84 & 0.00 & 165.91\\
H-All & 0.00 & 0.00 & 102.70 &  &  &  &  & 0.00 & 15.00 & 0.00 & 1151.38\\
\midrule
Synergies & 0.00 & 0.00 & 0.00 & 371.14 & 22357.57 & 0.00 & 0.00 & 0.00 & 5.44 & 0.00 & 0.74\\
\bottomrule
\end{tabular}
\end{table}

%% file: Tables/02_ResultsInstalledCapacities2040.tex
\begin{table}[h!]
\caption{Installed Capacities for all scenarios aggregated for the whole energy system design in 2030. Existing assets are not reported in this table.}
\label{tab:00_InstalledCapacities2040}
\resizebox{24cm}{!}{%
\begin{tabular}{l*{14}{S}}
\toprule
& \multicolumn{3}{c}{\ce{H2} pipeline} & \multicolumn{2}{c}{Electricity grid (AC)} & \multicolumn{2}{c}{Wind} & \multicolumn{1}{c}{PV} & \multicolumn{2}{c}{Battery} & \multicolumn{2}{c}{Electrolyser} & \multicolumn{1}{c}{Fuel Cell} & \multicolumn{1}{c}{\ce{H2} storage}\\
& \multicolumn{1}{c}{offshore} & \multicolumn{1}{c}{onshore (new)} & \multicolumn{1}{c}{onshore (repurposed)} & \multicolumn{1}{c}{AC} & \multicolumn{1}{c}{DC} & \multicolumn{1}{c}{offshore} & \multicolumn{1}{c}{onshore} & \multicolumn{1}{c}{} & \multicolumn{1}{c}{offshore} & \multicolumn{1}{c}{onshore} & \multicolumn{1}{c}{offshore} & \multicolumn{1}{c}{onshore} & \multicolumn{1}{c}{} & \multicolumn{1}{c}{}\\
& \multicolumn{1}{c}{(GWkm)} & \multicolumn{1}{c}{(GWkm)} & \multicolumn{1}{c}{(GWkm)} & \multicolumn{1}{c}{(GWkm)} & \multicolumn{1}{c}{(GWkm)} & \multicolumn{1}{c}{(GW)} & \multicolumn{1}{c}{(GW)} & \multicolumn{1}{c}{(GW)} & \multicolumn{1}{c}{(GWh)} & \multicolumn{1}{c}{(GWh)} & \multicolumn{1}{c}{(GW)} & \multicolumn{1}{c}{(GW)} & \multicolumn{1}{c}{(GW)} & \multicolumn{1}{c}{(GWh)}\\
\midrule
Reference &  &  &  &  & 19023.02 & 57.19 & 24.31 & 0.00 &  &  &  &  &  & \\
T-1 (only onshore) &  &  &  & 9370.64 & 20445.31 & 61.86 & 18.46 & 1.09 &  &  &  &  &  & \\
T-2 (only offshore) &  &  &  & 0.00 & 41886.01 & 40.59 & 29.89 & 0.05 &  &  &  &  &  & \\
T-3 (no border cross) &  &  &  & 7163.05 & 22955.13 & 61.62 & 15.85 & 3.16 &  &  &  &  &  & \\
T-All &  &  &  & 3462.74 & 38213.42 & 39.35 & 31.05 & 0.00 &  &  &  &  &  & \\
\midrule
S-1 (only onshore) &  &  &  &  & 19283.82 & 54.25 & 1.71 & 0.00 &  & 5.56 &  &  &  & \\
S-2 (only offshore) &  &  &  &  & 18513.78 & 54.17 & 1.69 & 0.00 & 5.68 &  &  &  &  & \\
S-All &  &  &  &  & 18793.32 & 54.17 & 1.7 & 0.00 & 0.00 & 5.54 &  &  &  & \\
\midrule
H-1 (only onshore) & 0.00 & 0.00 & 1944.59 &  & 21456.35 & 62.41 & 0.98 & 0.00 &  &  &  & 28.03 & 2.79 & 119.44\\
H-2 (only offshore) & 1337.55 & 0.00 & 1409.63 &  & 20294.91 & 56.45 & 1.61 & 0.00 &  &  & 9.81 &  & 2.86 & 64.36\\
H-3 (no storage) & 4556.39 & 0.00 & 2156.79 &  & 20984.89 & 64.19 & 2.36 & 0.00 &  &  & 0.34 & 29.53 & 2.69 & \\
H-4 (only local use) &  &  &  &  & 20686.22 & 62.33 & 1.24 & 0.00 &  &  &  & 26.50 & 2.80 & 287.39\\
H-All & 0.00 & 0.00 & 1921.25 &  & 21753.12 & 62.50 & 0.98 & 0.00 &  &  & 0.00 & 28.21 & 2.79 & 117.81\\
\midrule
Synergies & 0.00 & 0.00 & 0.00 & 2697.00 & 38976.89 & 39.78 & 31.46 & 0.00 & 0.00 & 0.00 & 0.00 & 11.97 & 0.00 & 9.95\\
\bottomrule
\end{tabular}
}
\end{table}